\documentclass[preprint,svgnames,dvipsnames,authoryear]{elsarticle}
\usepackage{dsfont}
\usepackage{amsmath,amsmath,amssymb,amsthm}
\usepackage{hyperref}
\usepackage{tabularx}
\usepackage[left=3cm,right=3cm,top=3cm,bottom=3cm]{geometry}
\usepackage[colorinlistoftodos]{todonotes}
\usepackage{setspace}
\usepackage[linesnumbered,ruled,vlined]{algorithm2e}
\usepackage{listings}
\usepackage[normalem]{ulem}

\makeatletter
\def\ps@pprintTitle{%
\let\@oddhead\@empty
\let\@evenhead\@empty
\def\@oddfoot{}%
\let\@evenfoot\@oddfoot}
\makeatother

\SetCommentSty{mycommfont}

\definecolor{backcolour}{rgb}{0.95,0.95,0.92}
\graphicspath{{img/}} 

\newtheorem{prop}{Property}

\newcommand{\LtwoNorm}[1]{
\left\lVert #1 \right\rVert_{L^2}}

\newcommand{\VarLR}[1]{
\mathbb{V} \left( #1 \right)}

\newcommand{\VarBigLR}[1]{
\mathbb{V} \Bigl( #1 \Bigr)}

\newcommand{\VarbigLR}[1]{
\mathbb{V} \bigl( #1 \bigr)}

\newcommand{\ie}{\textnormal{i.e.,}~}
\newcommand{\eg}{\textnormal{e.g.,}~}
\newcommand{\seg}{\textnormal{see, e.g.,}~}

\newcommand{\wrt}{\textnormal{w.r.t.}~}
\newcommand{\aka}{\textnormal{a.k.a.}~}
\newcommand{\iid}{\textnormal{i.i.d.}~}
\newcommand{\myeqref}[1]{Eq.~(\ref{#1})}
\newcommand{\myfigref}[1]{Figure~\ref{#1}}
\newcommand{\myalgref}[1]{Algorithm~\ref{#1}}

\newcommand{\rev}[1]{#1}
\renewcommand{\sout}[1]{\unskip}
\newcommand{\del}[1]{{\color{red} \sout{#1}}} 

\begin{document}

\onehalfspacing 

\begin{frontmatter}
\title{Developments and applications of Shapley effects to reliability-oriented sensitivity analysis with correlated inputs}

\author[a,b,c]{Marouane Il Idrissi}
\author[a,b]{Vincent Chabridon}
\author[a,b,c,d]{Bertrand Iooss}

\address[a]{EDF Lab Chatou, 6 Quai Watier, 78401 Chatou, France}
\address[b]{SINCLAIR AI Lab., Saclay, France}
\address[c]{Institut de Mathématiques de Toulouse, 31062 Toulouse, France}
\address[d]{Corresponding Author - Email: bertrand.iooss@edf.fr - Phone: +33130877969}
\begin{abstract}
Reliability-oriented sensitivity analysis methods have been developed for understanding the influence of model inputs relative to events which characterize the failure of a system (e.g., a threshold exceedance of the model output). In this field, the target sensitivity analysis focuses primarily on capturing the influence of the inputs on the occurrence of such a critical event. This paper proposes new target sensitivity indices, based on the Shapley values and called ``target Shapley effects'', allowing for interpretable sensitivity measures under dependent inputs. Two algorithms (one based on Monte Carlo sampling, and a given-data algorithm based on a nearest-neighbors procedure) are proposed for the estimation of these target Shapley effects based on the $\ell^2$ norm. Additionally, the behavior of these target Shapley effects are theoretically and empirically studied through various toy-cases. Finally, the application of these new indices in two real-world use-cases (a river flood model and a COVID-19 epidemiological model) is discussed.

\end{abstract}

\end{frontmatter}
\section{Introduction}

Nowadays, numerical models are extensively used in all industrial and scientific disciplines to describe physical phenomena (\eg systems of ordinary differential equations in ecosystem modeling, finite element models in structural mechanics, finite volume schemes in computational fluid dynamics) in order to design, analyze or optimize various processes and systems. \rev{These numerical models are often useful from either a scientific standpoint (\eg by improving the understanding of modeled physical phenomena) or from an engineering standpoint (\eg to better assist a decision-taking process).} In addition to this tremendous growth in computational modeling and simulation, the identification and treatment of the multiple sources of uncertainties has become an essential task from the early design stage to the whole system life cycle. As an example, such a task is crucial in the management of complex systems such as those encountered in energy exploration and production \citep{DeRocquigny_Devictor_book_2008} and in sustainable resource development \citep{bev08}.

In addition, the emergence of global sensitivity analysis (GSA) of model outputs played a fundamental role in the development and enhancement of these numerical models (\seg \cite{piabev16,razjak20} for recent reviews). Mathematically, if the model inputs (resp. output) are denoted by $X$ (resp. $Y$) and the model is written $G(\cdot)$, such as
\begin{equation}
\label{eq:num_model}
Y = G(X) ,
\end{equation}
GSA aims at understanding the behavior of $Y$ with respect to (\wrt) $X = (X_1, \dots, X_d)^{\top}$ the vector of $d$ inputs. GSA has been extensively used as a versatile tool to achieve various goals: for instance, quantifying the relative importance of inputs regarding their influence on the output (\aka "ranking"), identifying the most influential inputs among a large number of inputs (\aka screening) or analyzing the input-output code \rev{(\ie the numerically modeled phenomenon)} behavior \citep{Saltelli_THE_PRIMER_book,ioolem15}. 

When complex systems are critical or need to be highly safe, numerical models can also be of great help for risk and reliability assessment \citep{Lemaire2009}. Indeed, to track potential failures of a system (which could lead to dramatic environmental, human or financial consequences), numerical models allow a simulation of its behavior far from its nominal one (\seg \cite{ricvit19} in flood hazard assessment). In such a context, analytical or experimental approaches can be inappropriate, too expensive, or too difficult to perform. Based on numerical simulations, the tail behavior of the output distribution can be studied and typical \emph{risk measures} can be estimated \citep{Rockafellar_Royset_JRUES_2015}. Among others, the probability that the output $Y$ exceeds a given threshold value $t \in \mathbb{R}$, given by $\mathbb{P}(Y > t)$ and often called a \emph{failure probability}, is widely used in many applications. When $\{Y > t\}$ is a rare event (\ie associated to a very low failure probability), advanced sampling-based or approximation-based techniques \citep{MorioBalesdent2015} are required to accurately estimate the failure probability. In this very specific context, dedicated sensitivity analysis methods have been developed, especially in the structural reliability community (\seg \cite{Wu1994, Song2009, Wei_Lu_CPC_2012}). In such a framework, called \emph{reliability-oriented sensitivity analysis} (ROSA) \citep{chabridon_phd_2018, Perrin_Defaux_JSC_2019}, the idea is to provide importance measures dedicated to the problem of rare event estimation. 

Formally, standard GSA methods mostly focus on quantities of interest (QoI) characterizing the central part of the output distribution (\eg the variance for Sobol' indices \citep{Sobol_MMCE_1993}, the entire distribution for moment-independent indices \citep{borgonovo_new_2007}),  while ROSA methods focus on risk measures and their associated practical difficulties (\eg costly to estimate, inducing a conditioning on the distributions, non-trivial interpretation of the indices). Following \cite{raguet_target_2018}, ROSA methods can be categorized regarding the type of study they consider, \ie according to the following two categories:
\begin{itemize}
\item \emph{target sensitivity analysis} (TSA) aims at measuring the influence of the inputs (considering their entire input domain) on the \emph{occurrence} of the failure event. This means considering the following random variable, defined by the indicator function of the failure domain: $\mathds{1}_{\{G(X) > t\}}$;
\item \emph{conditional sensitivity analysis} aims at studying the influence of the inputs on the \emph{conditional} distribution of the output $Y|\{G(X) > t\}$, \ie exclusively within the critical domain. By \myeqref{eq:num_model}, a conditioning also appears on the inputs' domain.
\end{itemize}
Various indices have been proposed to tackle these two types of studies (\seg \cite{Li_Lu_SS_2012, Wei_Lu_CPC_2012, Perrin_Defaux_JSC_2019, marrel_statistical_2020}). The present paper is dedicated to ROSA (under the assumption that the QoI is a failure probability) and focuses on a TSA study. However, a new consideration for TSA is addressed in the present work: the possible statistical dependence between the inputs.

Indeed, most of the common GSA methods (and it is similar for the ROSA ones) have been developed under the assumption of independent inputs. As an example, the well-known Sobol' indices \citep{Sobol_MMCE_1993} which rely on the so-called functional analysis of variance (ANOVA) and Hoeffding decomposition \citep{Hoeffding_1948}, can be directly interpreted as shares of the output variance that are due to each input and combination of inputs (called ``interactions'') as long as the inputs are independent. 

When the inputs are dependent, the inputs' correlations dramatically alter the interpretation of the Sobol' indices. To handle this issue, several approaches have been investigated in the literature. For instance, \cite{jacques_sensitivity_2006} proposed to estimate indices for groups of correlated inputs. However, this approach does not allow for a quantification of the influence of individual inputs. Amongst other similar works, \cite{lirab10,chastaing_generalized_2012} proposed to extend the functional ANOVA decomposition to a more general one (\eg taking the covariance into account). However, the indices obtained for these approaches can be negative, which limits their practical use due to interpretability challenges \rev{(\ie as a share of the output's variance)}. In addition to this, other works (\seg \cite{xu_uncertainty_2008,martar12}) considered a Gram–Schmidt procedure to decorrelate the inputs and proposed to estimate two kinds of contributions for each variable (an uncorrelated one and a correlated one). These developments finally resulted in the proposition of a set of four Sobol' indices (instead of the two standard ones which are the first-order index and total index in the independent case) which enable the correlation effects to be fully captured in a GSA \citep{martar15}. Despite this achievement, this approach remains difficult to implement in practice (see \cite{beneli19} for extensive studies). 
Finally, the VARS approach \citep{doraz20} (allowing a thorough analysis of the inputs-output relationships) can handle input correlation but is out of scope of the present work which only focuses on variance-based sensitivity indices, directly computed from the numerical model.

Recently, another method has been developed by considering another type of indices: the \emph{Shapley effects}. The initial formulation originates from the ``Shapley values'' developed in the field of Game Theory \citep{Shapley_1953, osborne_course_1994}. The underlying idea is to fairly distribute both gains and costs to multiple players working cooperatively. By analogy with the GSA framework, the inputs can be seen as the players while the overall process can be seen as attributing shares of the output variability to the inputs. Considering the variance of the output in a GSA formulation leads to the so-called ``Shapley effects'' proposed by \cite{Owen_ASAJUQ_2014}. In the same vein, \citet{Owen_Prieur_ASAJUQ_2017, iooss_shapley_2019, beneli19} bridge the gap between Sobol' indices and Shapley effects while illustrating the usefulness of these new indices to handle correlated inputs in the GSA framework.

Thus, the present work attempts to to extend the use of Shapley effects to the ROSA context. Overall, the main objective is to develop a ROSA index which enables TSA to be performed (\ie capturing the influence of the inputs on a risk measure, typically a failure probability here) under the constraint of dependent inputs. This work relies on the use of recent promising results and numerical tools (both in field of TSA \cite{spagnol_phd_2020} and Shapley effects' estimation \cite{broto_variance_2020}).

The outline of this paper is the following. Section \ref{sec:gsa_dep} is devoted to a pedagogical introduction of the statistical dependence challenges for variance-based sensitivity indices, that can be solved by Shapley effects. Section \ref{sec:shapley_rosa} presents a new formulation of TSA, based on Shapley effects leading to the novel target Shapley effects, while Section \ref{sec:estimation} develops two algorithms for their estimation. Section \ref{sec:toycases} provides illustrations on simple toy-cases which give analytical expressions of the target Shapley effects, allowing deeper appreciation of their behavior. Section \ref{sec:appli} applies these new sensitivity indices to two use-cases: a simplified model of a river flood and an epidemiological model applied to the COVID-19 pandemic. Finally, Section \ref{sec:ccl} gives conclusions and research perspectives.

Throughout this paper, the mathematical notation $\mathbb{E}(\cdot)$ (resp. $\mathbb{V}(\cdot)$) will represent the expectation (resp. variance) operator.

\section{Variance-based sensitivity analysis with dependent inputs: the Shapley solution}
\label{sec:gsa_dep}

While devoted to computer experiments, GSA has close connections with multivariate data analysis and statistical learning \citep{chr90,hastib02}.
Indeed, in all these topics, one important issue is often to provide a weight to some variables (the inputs) \wrt its impact on another variables (the outputs). Depending on the domain, such a weight can either be called a ``sensitivity index'' or an ``importance measure''. A very convenient way is to base these weights on the ANOVA (analysis of variance) decomposition \citep{chr90,Sobol_MMCE_1993} of the output variance. Indeed, such a decomposition provides a natural division of the output's variance in shares attributed to each input. The principle of the ``variance-based sensitivity indices'' \citep{Saltelli_THE_PRIMER_book} consists then in understanding how to separate the contribution of each $X_i$ from the variance of $Y$. However, due to potential statistical dependencies between inputs, this decomposition cannot be directly performed. Starting from a simple example of a linear model, chosen for pedagogical purposes, this section provides a reminder on this topic while illustrating the important potential of Shapley effects in practice.

\subsection{Understanding the correlation issues via the linear model case}\label{sec:linear}

In this section, the aim is to quantify the relative importance of $d$ scalar inputs $X_j$ ($j=1,\ldots,d$) by fitting on a data sample (coming from the model \myeqref{eq:num_model}) a linear regression model so as to predict a scalar output $Y$:
\begin{equation}\label{eq:linmodel}
Y(X) = \sum_{j=0}^d \beta_j X_j + \epsilon \;,
\end{equation}
where $X_0=1$, $\mathbf{\beta}=(\beta_0,\ldots,\beta_d)^{\top} \in \mathbb{R}^{d+1}$ is the effects vector and $\epsilon \in \mathbb{R}$ the model's error of variance $\sigma^2$. If a sample of inputs and outputs $\displaystyle (\mathbf{X}^n,\mathbf{Y}^n) = \left(X_1^{(i)},\ldots,X_d^{(i)},Y^{(i)}\right)_{i=1,\dots,n}$ is available (with $n > d$), the Ordinary Least Squares method (see, e.g., \cite{chr90}) can easily be used to estimate the parameters $\beta$ and $\sigma^2$ in the linear regression model in \myeqref{eq:linmodel}. Moreover, one obtains the predictor $\widehat{Y}(x^*)$ of $Y$ at any prediction point $x^{*}$. An important validation metric of this model is the classical \emph{coefficient of determination} given by:
\begin{equation}
R^2_{Y(X)} = \sum_{i=1}^{n} \left. \left[ \widehat{Y}(X^{(i)}) - \bar{Y} \right]^2 \right/ \left[ Y^{(i)} - \bar{Y} \right]^2
\end{equation}
where $\bar{Y}$ is the output empirical mean. $R^2_{Y(X)}$ represents the percentage of output variability explained by the linear regression model of $Y$ on $X$. Finally, from \myeqref{eq:linmodel}, the variance decomposition expresses as:
\begin{equation}\label{eq:ancova}
\mathbb{V}(Y) = \sum_{j=1}^d \beta_j^2 \mathbb{V}(X_j) + 2 \sum_{k>j} \beta_j \beta_k \mathrm{Cov}(X_j,X_k) + \sigma^2 \;.
\end{equation}

In the specific case of independent inputs, the covariance terms cancel and the standard ANOVA (\ie $\mathbb{V}(Y) = \sum \beta_j^2 \mathbb{V}(X_j) +\sigma^2$) is obtained. Then, global sensitivity indices, called Standardized Regression Coefficients (SRC), can be directly computed:
\begin{equation}
\operatorname{SRC}_{j}=\beta_{j}\sqrt{\mathbb{V}(X_{j}) / \mathbb{V}(Y)}\;.
\label{eq:SRC}
\end{equation}
The estimation of the SRC is made by replacing the terms in \myeqref{eq:SRC} by their estimates. Interestingly, this metric for relative importance is signed (thanks to the regression coefficient sign), giving the sense of variation of the output \wrt each input. Moreover, $\operatorname{SRC}_{j}^2$ represents a share of variance and the sum of all the $\operatorname{SRC}_{j}^2$ approaches $R^2$ (\ie the amount of explained variance by the linear model). Note that, in a perfect linear regression model (\ie without any random error term $\epsilon$), $\operatorname{SRC}_{j}$ is equal to the linear Pearson's correlation coefficient between $X_j$ and $Y$ (denoted by $\rho(X_j,Y)$).
Note also that the ANOVA and SRC$^2$ extend to the functional ANOVA and Sobol' indices in the general (non-linear model) case (see \ref{app:sobol}).

When the inputs are dependent, the main concern is to allocate the covariance terms in \myeqref{eq:ancova} to the various inputs. In this case, the Partial Correlation Coefficient (PCC) has been promoted in GSA \citep{heljoh06,Saltelli_THE_PRIMER_book} as a substitute to the SRC, in order to cancel the effects of other inputs when allocating the weight of one input $X_j$ in the variance of $Y$:
\begin{equation}
\operatorname{PCC}_{j} = \rho(X_{j}-\widehat{X_{-j}},Y-\widehat{Y_{-j}})
\end{equation}
where $X_{-j}$ is the vector of all the $d$ inputs except $X_j$, $\widehat{X_{-j}}$ is the prediction of the linear model expressing $X_{j}$ \wrt $X_{-j}$ and $\widehat{Y_{-j}}$ is the prediction of the linear model $Y$ \wrt $X_{-j}$. 
However, PCC is not a right sensitivity index of the input. 
Indeed, it consists in measuring the linear correlation between $Y$ and $X_j$ by fixing $X_{-j}$, and  is then a measure of the linearity (and not the importance) between the output and one input.

Instead of controlling other inputs $X_{-j}$ such as done in the PCC, the Semi-Partial Correlation Coefficient (SPCC) quantifies the proportion of the output variance explained by $X_j$ after removing the information brought by $X_{-j}$ (on $X_j$) \citep{johleb04}:
\begin{equation}\label{eq:SPCC}
\operatorname{SPCC}_{j} = \rho(X_{j}-\widehat{X_{-j}},Y)  \;.
\end{equation}
SPCC can also be expressed by using the relation $\operatorname{SPCC}_j^2 = R^2_{Y(X)} - R^2_{Y(X_{-j})}$, which clearly shows that SPCC gives the additional explanatory power of the input $X_j$ in the linear regression model of $Y$ on $X$. However, the SPCC of highly correlated inputs will be small, despite their ``real'' explanatory power on the output. This aspect seems to be the main drawback of SPCC and probably explains its lack of popularity for GSA purposes.

\begin{figure}[!ht]
    \centering
    \includegraphics[scale=.8]{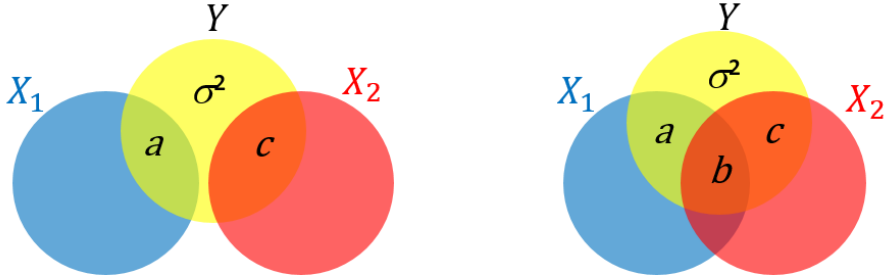}
    \caption{Inspired from \cite{clo19}. Illustration scheme of the effect of two inputs $X_1$ and $X_2$ on an output variable $Y$ when they are: uncorrelated (left) or correlated (right).}
    \label{fig:venn}
\end{figure}

\rev{To give an intuitive view of the limitations induced by the \emph{multicollinearity} of the inputs (\ie when inputs are linearly correlated to each other), Venn diagrams can be used (see, \myfigref{fig:venn}) in the case of two inputs, $X_1$ and $X_2$, and one output $Y$.} From \myfigref{fig:venn}, the coefficient of determination can be written as:
\begin{equation}
R^2_{Y(X_1,X_2)} = \frac{a+b+c}{a+b+c+\sigma^2} \;,
\end{equation}
where $a+b+c+\sigma^2$ is equal to the variance of $Y$ and $a+b+c$ represents the part of explained variance by the regression model (with $b=0$ in the uncorrelated case). In this elementary example, the previously introduced sensitivity indices are given by \cite{clo19}:
\begin{equation}
\begin{array}{lcllcl}
\mbox{SRC}_1^2 &=& (a+b)/(a+b+c+\sigma^2) \;, & \mbox{SRC}_2^2 &= &(c+b)/(a+b+c+\sigma^2) \;,\\
\mbox{PCC}_1^2&=&a/(a+\sigma^2) \;, & \mbox{PCC}_2^2 &=& c/(c+\sigma^2) \;,\\
\mbox{SPCC}_1^2 &=& a/(a+b+c+\sigma^2) \;, & \mbox{SPCC}_2^2 &=& c/(a+b+c+\sigma^2) \;.
\end{array}
\end{equation}
Thus, one can understand the limitations of SRC, PCC and SPCC when correlation is present: the variance share which comes from the correlation between inputs (\ie the $b$ value in \myfigref{fig:venn} - right) is allocated two times with the SRC but not allocated at all with SPCC, while PCC does not represent any variance sharing.

The three problems above can be solved by using another sensitivity index which finds a way to partition the $R^2$ among the $d$ inputs: the LMG \citep{linmer80,gro06} (acronym based on the authors' names, \ie ``Lindeman - Merenda - Gold'') uses sequential sums of squares from the linear model and obtains an overall measure by averaging over all orderings of inputs. Mathematically, let \rev{$A$} be a subset of indices in the set of all subsets of $\{1,\ldots,d\}$ and \rev{$X_A=(X_j:j \in A)$} a group of inputs. LMG is based on the measure of the elementary contribution of any given variable $X_j$ to a given subset model \rev{$Y(X_A)$} by the increase in $R^2$ that results from adding that predictive variable to the regression model:
\begin{equation}
\operatorname{LMG}_{j} = \frac{1}{d!} \sum_{\substack{\pi \in \mbox{\scriptsize permutations}\\ \mbox{\scriptsize of } \{1,\ldots,d\}}} \left[ R^2_{Y(X_{v \cup \{j\}})} - R^2_{Y(X_v)} \right]
\label{eq:LMGperm}
\end{equation}
with $v$ the indices entered before $j$ in the order $\pi$. 
In \myeqref{eq:LMGperm}, the sum is performed over all the permutations of $\{1,\ldots,d\}$.
For the case of two inputs (see \myfigref{fig:venn}), we can easily show that:
\begin{equation}
\operatorname{LMG}_1 = (a+b/2)/(a+b+c+\sigma^2) \;, \; \operatorname{LMG}_2 = (c+b/2)/(a+b+c+\sigma^2) \;.
\end{equation}
Then, in the LMG framework, the $R^2_{Y(X_1,X_2)}$ has been perfectly shared into two parts with an equitable distribution of the $b$ term between $X_1$ and $X_2$.


This allocation principle exactly corresponds to the application of the \emph{Shapley values} \citep{Shapley_1953} on the linear model. 
This attribution method has been primarily used in cooperative game theory, allowing for a cooperative allocation of resources between players based on their collective production (see \ref{app:def_shap_val} for a more formal definition).
The Shapley values solution consists in fairly distributing both gains and costs to several actors working in coalition.
In situations when the contributions of each actor are unequal, it ensures that each actor gains as much or more as they would have from acting independently. 
Now, if the actors are identified with a set of inputs and the value assigned to each coalition is identified to the explanatory power of the subset of model inputs composing the coalition, one obtains the LMG in \myeqref{eq:LMGperm}.

\subsection{Shapley effects}
\label{subsec:shapley_eff}

In the general case, \rev{when no assumption is made} on the model $G(\cdot)$ (see \myeqref{eq:num_model}), variance-based sensitivity indices have been developed (see, \cite{Sobol_MMCE_1993,Saltelli_THE_PRIMER_book}) and applied to perform a GSA of complex models (\seg \cite{nosels11}).

\rev{When the inputs are assumed to be independent, they} allow the variance of the model output to be decomposed \rev{according to each possible subsets of inputs} (called ``Sobol' indices''). \rev{They are a means to measure the individual effects of inputs, as well as the effect of their interaction} (see, \ref{app:sobol} for the theoretical details).

\rev{When the inputs are effectively dependent, the Sobol' indices lose their inherent interpretation (\ie decomposition in individual and interaction effects). To remedy this drawback, \citet{Owen_ASAJUQ_2014} recently proposed game theoretic GSA indices, in the same fashion as the LMG indices (see, \myeqref{eq:LMGperm}), inspired by the Shapley values of cooperative games. They are defined by:}
\begin{equation}
\label{eq:shapley}
Sh_j = \frac{1}{d}  \sum_{A \subset \{-j\} } {d-1 \choose |A|}^{-1} \left( \textrm{val}(A \cup \{j\} ) - \textrm{val}(A) \right),
\end{equation}
where $\textrm{val}(A)$ is called \emph{cost} function \rev{(or \emph{value} function)} assigned to a subset $A \in \mathcal{P}_d$ of inputs, $\mathcal{P}_d$ \rev{denotes} the set of all possible subsets of $\{1, \dots, d\}$, $\{-j\}$ denotes the set of indices $\{1, \dots, d\} \setminus j$ and $|A|$ \rev{denotes} the cardinal \rev{number} of $A$.

For GSA purposes, \cite{Owen_ASAJUQ_2014} proposes to use the ``closed Sobol' indices'' as the value function in \myeqref{eq:shapley}:
\begin{equation}
\label{eq:cost_fun}
\textrm{val}(A) = S_A^{\textrm{clos}} = \frac{\mathbb{V}\left(\mathbb{E}\bigl[ G(X) \bigm| X_A\bigr]\right)}{\mathbb{V}\left(G(X)\right)},
\end{equation}
The attribution properties of the Shapley values applied to this particular cost function, leads to the definition of the \emph{Shapley effects}:
\begin{equation}
Sh_j = \frac{1}{d}  \sum_{A \subset \{-j\} } {d-1 \choose |A|}^{-1} \left( S^{\textrm{clos}}_{A \cup \{j\}} - S^{\textrm{clos}}_A \right).
\label{eq:shap_eff}
\end{equation}
These indices allow for a quantification of \rev{the importance of each input}, which intrinsically takes into account both interaction and dependence \rev{effects on the output's variance}. Moreover, two important properties of the Shapley effects allow \rev{for their interpretation}: they sum up to one and are non-negative. \rev{Thus, they can be considered as a decomposition of the output's variance}. They allow for input ranking by \rev{attributing to} each input a percentage of the \rev{variable of interest's} variance. These indices have been extensively studied \rev{by} \cite{Owen_Prieur_ASAJUQ_2017, iooss_shapley_2019}. An \rev{alternative} way of defining the Shapley effects has \rev{also} been proposed, by taking the following cost function:
\begin{equation}
\label{eq:cost_fun2}
\textrm{val}(A) =  \frac{\mathbb{E}\left[\VarbigLR{ G(X)| X_{\overline{A}} }\right]}{\mathbb{V}\left(G(X)\right)}
\end{equation}
where $\overline{A} = \{ 1,\dots, d \} \setminus A$. \rev{This alternative definition leads} to an equivalent definition of the Shapley effects \rev{(\myeqref{eq:shap_eff})}, as outlined by \rev{\cite{song_shapley_2016}, and allows for additional estimation methods.}

In order to illustrate Shapley effects' \rev{attributions}, one can first consider a model with three inputs $X=(X_1, X_2, X_3)^{\top}$. From \myeqref{eq:shap_eff}, one has:
\begin{align*}
Sh_1 &= \frac{1}{3}S_1^{\textrm{clos}} \\
&+ \frac{1}{6} \Bigl[ \bigl( S_{\{1,2 \}}^{\textrm{clos}} - S_2^{\textrm{clos}}\bigr) + \bigl( S_{\{1,3 \}}^{\textrm{clos}} - S_3^{\textrm{clos}}\bigr)\Bigr] \\
&+ \frac{1}{3}  \bigl(S_{\{1,2,3\}}^{\textrm{clos}} - S_{\{2,3 \}}^{\textrm{clos}} \bigr). 
\end{align*}
\rev{If the three inputs are assumed to be independent, this result leads to:}
\begin{equation*}
\label{eq:sh_indep}
Sh_1 = S_1 + \frac{1}{2}S_{\{1,2\}} + \frac{1}{2}  S_{\{1,3\}} + \frac{1}{3} S_{\{1,2,3\}}
\end{equation*}
where one can notice that the Shapley \rev{effects' decomposition} consists \rev{in allocating} the initial Sobol' index, plus an equal share of the interaction effects between all the inputs. However, if dependence between inputs is assumed, this behavior cannot be clearly illustrated, except \rev{when a linear model $G(\cdot)$ is assumed} (see Subsection \ref{sec:linear}).

\rev{The quantity} $\bigl( S^{\textrm{clos}}_{A \cup \{ j\}} - S^{\textrm{clos}}_A \bigr)$ can be interpreted as being a quantification of the \emph{\rev{marginal} effects of the input $j$ in relation to the subset of variables $A$}. \rev{It is heavily linked to the notion of \emph{marginal contributions} of cooperative games, aiming at quantifying the bargaining power of a player in an allocation process \citep{brandenburger_cooperative_2007}.} If $S^{\textrm{clos}}_A$ is believed to contain the initial effects of \rev{the inputs in} $A$, plus their interaction effects, and any effect due to \rev{their dependence structure}, then the increment \rev{$\bigl( S^{\textrm{clos}}_{A \cup \{ j\}} - S^{\textrm{clos}}_A \bigr)$} quantifies the initial effect of the input $j$, its interaction effects with \rev{the inputs in} $A$, and the effects due to their dependence. 
Then, the Shapley attribution weighs all \rev{the marginal} effects, in order to \rev{assess the effective influence of the involved inputs through their marginal contributions}, in the same fashion as the LMG \rev{indices for a} linear model, as depicted in Section \ref{sec:linear}.

\rev{It is important to note that the above mentioned alleged decomposition of $S_A^{\textrm{clos}}$ cannot be verified due to the lack of a univocal functional variance decomposition when inputs are dependent. However, the interpretation of the Shapley effects does not rely on the chosen cost function to be meaningful, but rather on its ability to quantify marginal contributions. Empirical studies and analytical studies show that the choice of $S_A^{\textrm{clos}}$ as a cost function remains pertinent, even when inputs are dependent (see, \cite{iooss_shapley_2019}).}
\section{Reliability-oriented Shapley effects for target sensitivity analysis}
\label{sec:shapley_rosa}

\subsection{A brief overview of reliability-oriented sensitivity analysis}
\label{subsec:tsa}

When focusing on complex systems, one often needs to prepare for possible critical events, which potentially have a low occurence probability but \rev{lead to a} system \rev{failure}. Such failures may have dramatic human, environmental and economic consequences, depending on the context. \rev{The fields} of reliability assessment and risk analysis \citep{Lemaire2009,ricvit19}, \rev{aim to prevent these failures}. Mathematically, \rev{a reliability} problem focuses on a \emph{risk measure} computed from the tail of the variable of interest's distribution \citep{Rockafellar_Royset_JRUES_2015}. Performing sensitivity analysis in such a context requires \rev{the} use \rev{of} dedicated tools, which have been \rev{developed} by various authors under the denomination of ``reliability-oriented sensitivity analysis'' (ROSA) (\seg \cite{Perrin_Defaux_JSC_2019,Derennes_Morio_Simatos_MACS_2021,marrel_statistical_2020}). A large panel of ROSA methods have been proposed in the structural reliability community such as, for example, several variance-based approaches (\seg \cite{Morio_SMPT_2012,Wei_Lu_CPC_2012,Perrin_Defaux_JSC_2019,Chabridon_Chapter_2020}) and moment-independent approaches (\seg \cite{Cui_Lu_SCTS_2010,Li_Lu_SS_2012,Derennes_Morio_Simatos_MACS_2021}). From the GSA community, several extensions have also been proposed in order to study risks, or reliability measures. The contrast-based indices proposed by \cite{FortKleinRachdi2016} are, amongst others, an example of a versatile tool which can handle several types of QoI. They were applied in the works of \cite{BrowneFortIoossLeGratiet_sensi2017,MaumeDeschampsNiang_2018} in quantile-oriented formulations. Other formulations such as the quantile-based global sensitivity measures \citep{Kucherenko_Song_Wang_RESS_2019} or other indices related to dependence measures \citep{raguet_target_2018,marrel_statistical_2020} have been proposed.

In the context of reliability assessment, a typical risk measure is the \emph{failure probability} given by:
\begin{equation}
p_t^Y \overset{\textrm{def}}{=} \mathbb{P} \left( Y > t \right) = \mathbb{P} \left( G(X) > t \right) = \mathbb{E} \left[ \mathds{1}_{\{G(X)>t\}}(X) \right] = \mathbb{E} \left[ \mathds{1}_{\mathcal{F}_t}(X) \right]
\label{eq:tsa_qoi}
\end{equation}
where \rev{$t \in \mathbb{R}$} represents a threshold characterizing the state of the system. Typically, the event $\{Y > t\}$ denotes \rev{a} failure event \rev{(\ie the system described by the model $G(X)$ enters a failure state)}. As for $\mathcal{F}_t$, it represents the input failure domain, \ie $\mathcal{F}_t \overset{\textrm{def}}{=}\{X~|~G(X) > t\}$.

Performing a ROSA study poses a few challenges: firstly, the variable of interest here is \rev{not directly} $Y$ \rev{anymore}, but \rev{rather} a binary \rev{random variable} whose occurrence is characterized by the indicator function $\mathds{1}_{\mathcal{F}_t}(X)$; secondly, \rev{in practice, these failure events are typically} ``rare events'', associated to a low failure probability which might be difficult to estimate in practice \rev{through typical sampling methods} \citep{MorioBalesdent2015}; thirdly, the type of study one desires to perform has to be reinterpreted regarding the new QoI. Regarding this last point, \cite{raguet_target_2018} focus on two types of studies when dealing with critical events: the first one, called \emph{target sensitivity analysis} (TSA), aims at catching the influence of the inputs on the occurrence of the failure event, while the second one, called ``conditional sensitivity analysis'' aims at studying the influence of the inputs once the threshold value has been reached (\ie within the failure domain). The present paper is dedicated to ROSA (under the assumption that the QoI is a failure probability given by \myeqref{eq:tsa_qoi}) and \rev{aims at developing tools for} TSA.

To illustrate this \rev{paradigm }in plain text, one can refer to the study of the water level in a river protected by a dyke. From the traditional GSA point of view, the central question would be ``Which inputs influence the water level?'', while in the TSA paradigm, one focuses more on the question ``Which inputs influence the occurrence of a flood?''. Note that this particular example is studied \rev{in depth} in Subsection~\ref{subsec:flood_case}.

\begin{figure}[!ht]
    \centering
    \includegraphics[width=\textwidth]{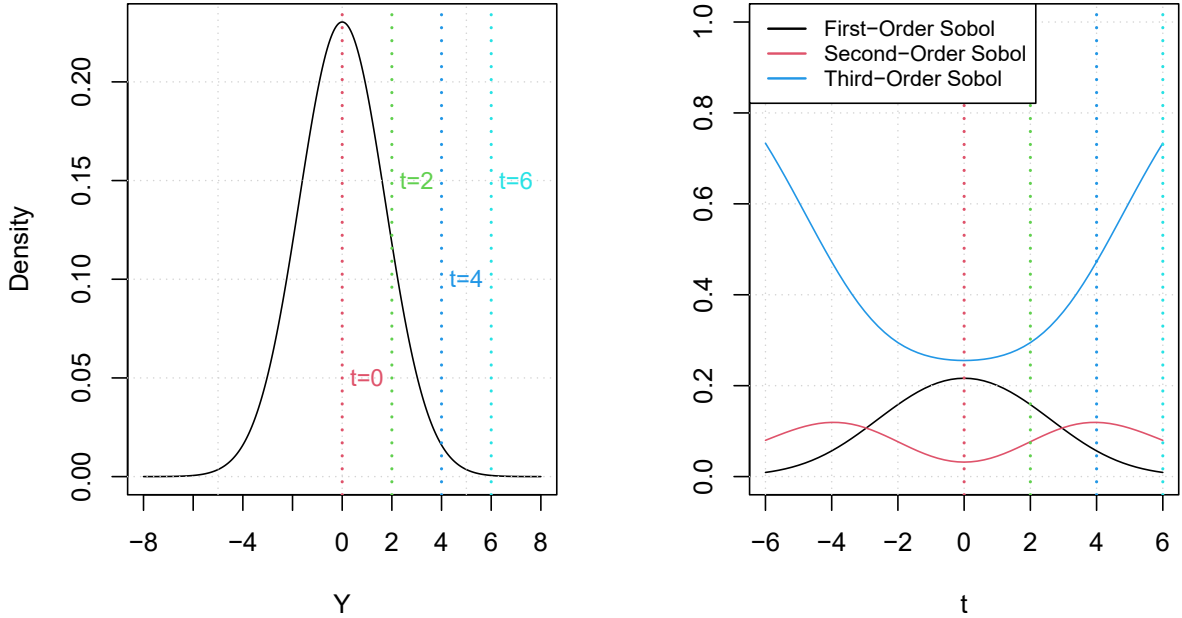}
    \caption{Probability density function of the output with four different threshold values (left) and the related target Sobol' indices (right) for $Y$ being the sum of three independent Gaussian random variables.}
    \label{fig:tsa_sobolex}
\end{figure}

\rev{When the inputs are assumed to be independent}, a first category of sensitivity indices dedicated to TSA are the ``target Sobol' indices'' whose first formulation has been proposed by \cite{Li_Lu_SS_2012}. \rev{In this document, the original Sobol' indices, in the TSA context, are denoted:
\begin{equation}
\textrm{T-S}_A = \sum_{B \subseteq A} (-1)^{|A| - |B|}  \frac{\mathbb{V}\left(\mathbb{E}\left[ \mathds{1}_{\mathcal{F}_t}(X)~|~X_B\right]\right)}{\mathbb{V}\left(\mathds{1}_{\mathcal{F}_t}(X)\right)}
\label{eq:original_tsasobol}
\end{equation}
where $\mathbb{V}\left(\mathds{1}_{\mathcal{F}_t}(X)\right) = p_t^Y(1- p_t^Y)$. Similarly, the closed Sobol' indices (see \ref{app:sobol}) are defined as follows:
\begin{equation}
\textrm{T-S}_A^{\textrm{clos}} = \frac{\mathbb{V}\left(\mathbb{E}\left[ \mathds{1}_{\mathcal{F}_t}(X)~|~X_A\right]\right)}{\mathbb{V}\left(\mathds{1}_{\mathcal{F}_t}(X)\right)}.
\label{eq:tsa_sobol}
\end{equation}
}
Several estimation schemes for these indices have been proposed \rev{when dealing with rare failure events} \citep{Wei_Lu_CPC_2012,Perrin_Defaux_JSC_2019}. To illustrate the behavior of \rev{these indices}, one can consider \rev{a linear} model given by $Y=X_1 + X_2 + X_3$, with $X=(X_1, X_2, X_3)^{\top}$, three standard Gaussian random variables assumed to be independent. The left plot of \myfigref{fig:tsa_sobolex} represents the probability density function (pdf) of $Y$, \rev{along with four different threshold values}, corresponding to four different failure probability levels. The right plot of \myfigref{fig:tsa_sobolex} presents the different values of $\textrm{T-S}_A$, \wrt the threshold $t$. Note that, \rev{when the inputs are assumed to be independent}, the second-order \rev{Sobol'} indices \rev{(\ie when $|A|=2$)} verify $\textrm{T-S}_{\{1,2\}}=\textrm{T-S}_{\{1,3\}}=\textrm{T-S}_{\{2,3\}}$. \rev{One can additionally remark} that \rev{as soon as $t$ induces a low or high failure probability} (\ie ``close'' to $0$ or $1$), the third-order \rev{(\ie $|A|=3$)} closed Sobol' index for TSA increases, indicating high interaction between the three inputs. Note that \rev{studying this behavior} falls under the conditional sensitivity analysis paradigm, which is out of the scope of this paper. However, the acknowledgment of this phenomenon remains important for better understanding \rev{the proposed indices' behavior}.

Another category of sensitivity indices dedicated to TSA \rev{fall under the category of} ``moment-independent'' ROSA indices. Among others, one can mention the two indices proposed by \cite{Cui_Lu_SCTS_2010} which are given by:
\begin{subequations}
\begin{alignat}{3}
	\eta_A &= \frac{1}{2} \mathbb{E} \left[ \left| p_t^Y - p_t^{Y|X_A} \right| \right]\\
	\delta_A &= \frac{1}{2} \mathbb{E} \left[ \left( p_t^Y - p_t^{Y|X_A} \right)^2 \right]
\end{alignat}
\label{eq:eta}
\end{subequations}
where $p_t^{Y|X_A}$ denotes the conditional failure probability when $X_A$ is fixed. Note that, if \rev{$\eta_A$} does not require any independence assumption for \rev{a meaningful quantification of} input influence, it is known to be difficult to estimate in practice \citep{Derennes_Morio_Simatos_MACS_2021}. As for \rev{$\delta_A$}, it is simply proportional to the target \rev{closed} Sobol' index given in \myeqref{eq:tsa_sobol}. Note that an extension of $\delta_A$ has been proposed in \cite{Li_Lu_AST_2016} for correlated inputs. It relies on a similar orthogonalization procedure strategy as proposed by \cite{martar12} for usual Sobol' indices. However, as mentioned previously, this tends to increase the number of estimated indices to properly interpret the \rev{inputs' influence}.

The following section \rev{aims at introducing} the \emph{distance-based TSA indices}, \rev{while }highlighting their links with existing TSA indices. \rev{New} TSA indices inspired from the Shapley \rev{values} \rev{(see, Section \ref{subsec:shapley_eff}) are then proposed.}

\subsection{Distance-based TSA indices}

As \rev{outlined} by several authors \citep{FortKleinRachdi2016,raguet_target_2018}, it can be noted that $\mathbb{V}\left(\mathbb{E}\left[ Y|X_A\right]\right) = \mathbb{E}\left[ \left( \mathbb{E}\left[ Y|X_A\right] -  \mathbb{E}\left[Y\right]\right)^2\right]$. This equality can be interpreted as the expected squared distance between two expectations, and thus allows to apprehend \rev{closed} Sobol' indices (see \myeqref{eq:cost_fun}) as \rev{a particular case of} distance-based indices. \rev{This broader} point of view has been adopted by \cite{FortKleinRachdi2016} to provide a generalization of the Sobol' indices using \emph{contrast functions}.

By applying a similar idea \rev{for TSA}, one can extend the standard $\textrm{T-S}_A$ and $\eta_A$ \rev{definition} to more general cases \rev{based on} distances. One can \rev{then} define \rev{more general} distance-based TSA \rev{indices}, relative to a subset of inputs $A \in \mathcal{P}_d$, as follows:
\begin{equation}
\label{eq:D_dist}
\textrm{T-S}_A^{\mathcal{D}} = \frac{\mathbb{E}\Bigl[\mathcal{D}\left(p_t^Y, p_t^{Y|X_A}\right) \Bigr]}{\mathbb{E} \Bigl[ \mathcal{D}\left(p_t^Y, p_t^{Y|X}\right)\Bigr]}
\end{equation}
where $\mathcal{D}(\cdot,\cdot)$ \rev{can be} any distance function. Links can be made between \rev{this definition} and the indices presented previously, through the use of specific distance functions. For example, by choosing the distance derived from the $\ell^1$ norm (\ie the absolute difference), one can remark that the corresponding distance-based TSA index is proportional to the $\eta_A$ index:
\begin{equation}
\label{eq:l1_dist}
\textrm{T-S}^{\ell^1}_A = \frac{\mathbb{E} \Bigl[ \bigl| \mathbb{E}[\mathds{1}_{{F}_t} (X)] - \mathbb{E}[\mathds{1}_{{F}_t} (X) |X_A]\bigr|\Bigr]}{\mathbb{E}\Bigl[ \bigl\lvert \mathds{1}_{\mathcal{F}_t}(X) - \mathbb{E}\left[\mathds{1}_{\mathcal{F}_t}(X)\right] \bigr| \Bigr]}= \frac{2}{\mathbb{E}\Bigl[ \bigl\lvert \mathds{1}_{\mathcal{F}_t}(X) - \mathbb{E}\left[\mathds{1}_{\mathcal{F}_t}(X)\right] \bigr| \Bigr]} \eta_A  
\end{equation}

Moreover, by using the distance derived from the $\ell^2$ norm (\ie the squared difference), one can remark \rev{the resulting} distance-based TSA indices are equal to the closed Sobol' index for TSA, as defined in \myeqref{eq:tsa_sobol}:
\rev{
\begin{equation}
\label{eq:l2_dist}
\textrm{T-S}^{\ell^2}_A =\textrm{T-S}_A^{\textrm{clos}} = \frac{\VarBigLR{\mathbb{E}\bigl[ \mathds{1}_{\mathcal{F}_t}(X) \bigm| X_A\bigr]}}{\VarBigLR{\mathds{1}_{\mathcal{F}_t}(X)}}.
\end{equation}}
\rev{As outlined, the distance-based TSA indices are intimately related to existing ones (\ie $\eta_A$ and the closed Sobol' indices for TSA), and can thus be seen as a broader class of indices. Moreover, they are relevant candidates as cost functions for defining Shapley values inspired TSA indices, following a similar line of thinking as in Section \ref{subsec:shapley_eff}.}

\subsection{$(\mathcal{D})$-target Shapley effects}
In this subsection, a novel family of TSA indices is proposed, namely the $(\mathcal{D})$-target Shapley effects. As \rev{briefly} mentioned \rev{previously}, these indices are \rev{constructed} by taking distance-based TSA indices, defined in \myeqref{eq:D_dist}, as cost functions in a Shapley attribution procedure (see \myeqref{eq:shapley}). For a specific input $j \in \{1, \dots, d \}$, its $(\mathcal{D})$-target Shapley effects can be defined as being:
\begin{equation}
\label{eq:dtse_def}
\textrm{T-Sh}_j^{\mathcal{D}} = \frac{1}{d} \sum_{A \subset \{-j\}} {d-1 \choose |A|}^{-1} \left( \textrm{T-S}^{\mathcal{D}}_{A\cup \{j\}} - \textrm{T-S}^{\mathcal{D}}_A\right)
\end{equation}
where $\{ -j\} = \{ 1, \dots, d \} \setminus j$. The main property allowing for a clear interpretation of the $(\mathcal{D})$-target Shapley effects is the following:
\begin{prop}[$(\mathcal{D})$-target Shapley effects decomposition]
\label{prop:dtse_sum1}
Let $A \in \mathcal{P}_d$, and $val(A) = \textrm{T-S}_A^{\mathcal{D}}$. For any distance function $\mathcal{D}(.,.)$, the following property holds:
\begin{equation}
\label{eq:dtse_prop1}
\sum_{j=1}^d \textrm{T-Sh}_j^{\mathcal{D}} = 1.
\end{equation}
\end{prop}
It is important to note that this decomposition property does not rely on any independence assumption about the \rev{probabilistic model} of the inputs. However, in order to ensure a meaningful interpretation of these indices, \rev{(\ie as a percentage of a statistical dispersion)}, one needs to ensure that the $\textrm{T-Sh}_j$ are non-negative, for all $j=1,\dots, d$.

By choosing \rev{$\textrm{T-S}^{\ell^1}_A$ as a cost function (\ie $\mathcal{D}(x,y) = |x-y|$)}, one can then define the $(\ell^1)$-target Shapley effect associated to a variable $j \in \{1,\dots,d \}$ as being:
\begin{equation}
\label{eq:l1tse}
\textrm{T-Sh}_j^{\ell^1}=\frac{1}{d} \sum_{A \subset \{-j\}} {d-1 \choose |A|}^{-1} \left( \textrm{T-S}^{\ell^1}_{A \cup \{j\}} - \textrm{T-S}^{\ell^1}_{A} \right).
\end{equation}
These indices are non-negative (see the proof in \ref{ap:pos_proofl1}) which allows the $(\ell^1)$-target Shapley effects to be interpreted as the percentage of the mean absolute deviation of the indicator function (\ie $\mathbb{E}\Bigl[ \bigl\lvert \mathds{1}_{\mathcal{F}_t}(X) - \mathbb{E}\left[\mathds{1}_{\mathcal{F}_t}(X)\right] \bigr| \Bigr]$), \rev{allocated} to each input $X_j$, $j \in \{1, \dots,d\}$.

By choosing $\textrm{T-S}^{\ell^2}_A$ as a cost function (\ie $\mathcal{D}(x,y)=(x-y)^2$), the $(\ell^2)$-target Shapley effect associated to the variable $j \in \{1, \dots, d \}$ can be defined as:
\begin{align}
\label{eq:l2tse}
\textrm{T-Sh}_j^{\ell^2}&=\frac{1}{d} \sum_{A \subset \{-j\}} {d-1 \choose |A|}^{-1} \left( \textrm{T-S}^{\ell^2}_{A \cup \{j\}} - \textrm{T-S}^{\ell^2}_{A} \right).
\end{align}
Being also non-negative (see the proof in \ref{ap:pos_proofl2}), they can be interpreted as a percentage of the variance of the indicator function allocated to the input $X_j$, $j \in \{1, \dots,d\}$.
Moreover, using \myeqref{eq:l2_dist}, by analogy with Eqs. (\ref{eq:cost_fun}) and (\ref{eq:cost_fun2}) (in a similar fashion as the alternate cost function proposed by \cite{song_shapley_2016}), if one chooses to define the cost function $\textrm{val}(A)$ as being:
\begin{equation}
\label{eq:alt_costfun}
\textrm{T-E}_A \overset{\textrm{def}}{=} \frac{\mathbb{E}\Bigl[ \VarbigLR{\mathds{1}_{\mathcal{F}_t}(X) | X_{\overline{A}}} \Bigr]}{\VarBigLR{\mathds{1}_{\mathcal{F}_t}(X)}}
\end{equation}
with $\overline{A} = \{1, \dots, d\} \setminus A$, then one has an equivalent way of defining the $(\ell^2)$-target Shapley effect.

In the following, the $(\ell^2)$-target Shapley effect $\textrm{T-Sh}_j^{\ell^2}$ will be \rev{referred to as} ``the'' target Shapley effect and denoted $\textrm{T-Sh}_j$:
\begin{equation}
\textrm{T-Sh}_j \overset{\textrm{def}}{=} \textrm{T-Sh}_j^{\ell^2}.
\end{equation}

\section{Estimation methods and practical implementation of target Shapley effects}
\label{sec:estimation}

The estimation of the target Shapley effects \myeqref{eq:l2tse} can be \rev{done} into two \rev{distinct} steps:
\begin{itemize}
    \item \textbf{Step \#1}: estimation of the \emph{conditional elements}, \ie the estimation of \rev{$\textrm{T-S}^{\ell^2}_A$} or $\textrm{T-E}_A$ for all $A \in \mathcal{P}_d$;
    \item \textbf{Step \#2}: an \emph{aggregation procedure}, \ie a step to compute the $\textrm{T-Sh}_j$ by plugging in the previous estimations of Step \#1 in \myeqref{eq:l2tse}.
\end{itemize}
In the following, two estimation methods are proposed: the first one based on a Monte Carlo sampling procedure, and the second one based on a nearest-neighbor approximation technique.

\subsection{Monte Carlo sampling-based estimation}
\label{sec:mc_estim}

This procedure, introduced in \cite{song_shapley_2016} for the estimation of Shapley effects, relies on a Monte Carlo estimation of the conditional elements. It requires the ability to sample from the marginal distributions of the inputs (\ie $P_{X_A}$ for all $A \subseteq \{ 1, \dots, d\} \setminus \emptyset$), as well as from all the conditional distributions (\ie $P_{X_{\overline{A}}|X_A}$, for all possible subsets of inputs $A$). Additionally, one also needs to be able to evaluate the model $G(\cdot)$ which is \rev{usually} the case in the context of uncertainty quantification of numerical computer models (ignoring the potential difficulties related to the cost of a single evaluation of $G(\cdot)$) \citep{DeRocquigny_Devictor_book_2008}.

In order to estimate a conditional element \rev{$\textrm{T-S}_A^{\ell^2}$}, one needs to \rev{randomly} draw several \iid samples:
\begin{itemize}
    \item an \iid sample of size $N$ drawn from $P_X$ and denoted by $(X^{(1)}, \dots ,X^{(N)})$;
    \item another \iid sample of size $N_v$ drawn from $P_{X_A}$ and denoted by $(X_A^{(1)}, \dots, X_A^{(N_v)})$;
    \item for each element $X_A^{(i)}, i=1, \dots, N_v$, a corresponding sample of size $N_p$ drawn from $P_{X_{\overline{A}}|X_A}$ given that $X_A=X_A^{(i)}$ and denoted by $(\widetilde{X}_i^{(1)}, \dots, \widetilde{X}_i^{(N_p)})$.
\end{itemize}
Then, the Monte Carlo estimator of \rev{$\textrm{T-S}_A^{\ell^2}$} can be defined as:
\begin{equation}
\label{eq:l2_mc}
\widehat{\textrm{T-S}}_{A, \textrm{MC}} =  \frac{\sum_{i=1}^{N_v} \left( \frac{1}{N_p} \sum_{j=1}^{N_p} \mathds{1}_{\mathcal{F}_t}(\widetilde{X}_i^{(j)},X_A^{(i)})  - \widehat{p}_t^Y  \right)^2}{(N_v-1)\widehat{p}_t^Y(1-\widehat{p}_t^Y)}
\end{equation}
with
\begin{equation}
\label{eq:defprob_mc}
\widehat{p}_t^Y = \frac{1}{N} \sum_{i=1}^N \mathds{1}_{\mathcal{F}_t}(X^{(i)}).
\end{equation}
Finally, the aggregation procedure gives:
\begin{equation}
\widehat{\textrm{T-Sh}}_{j, \text{MC}} = \frac{1}{d} \sum_{A \subset \{-j\}} {d-1 \choose |A|}^{-1} \left( \widehat{\textrm{T-S}}_{A \cup \{ j\}, \text{MC}} - \widehat{\textrm{T-S}}_{A, \textrm{MC}} \right).
\label{eq:l2tse_mc}
\end{equation}
Thus, one gets that $\widehat{\textrm{T-Sh}}_{j, \text{MC}}$ is an unbiased consistent estimator of $\textrm{T-Sh}_j$. 

\myalgref{alg:mc_est} provides a detailed description on how to implement this estimator in practice. This estimation method requires $(N + d!\times (d-1)\times N_v \times N_p)$ calls to the numerical model $G(\cdot)$. \rev{Its empirical convergence \wrt $N_v$ is illustrated in \ref{app:empconv_mc}.}
As expected, this first estimation method can become quite expensive in practice. Moreover, numerical models usually encountered in industrial studies can be costly-to-evaluate, which can \rev{strongly limit the use of such a method in practice}.

\begin{algorithm}[!ht]
\SetKwFunction{sJ}{simJoint}
\SetKwFunction{sM}{simMarginal}
\SetKwFunction{sC}{simConditional}

\KwIn{$G, t, d, N, N_v, N_p,\sJ, \sM, \sC$}
\KwOut{$(\widehat{\textrm{T-Sh}}_{j, \textrm{MC}})_{j=1,\dots,d}$}
\tcc{Sample from the joint distribution}
$(X^{(1)}, \dots, X^{(N)})\leftarrow \texttt{sim\_joint}(N)$\\
\tcc{Estimate the failure probability}
$\widehat{p}_t^Y \leftarrow \frac{1}{N} \sum_{i=1}^N \mathds{1}_{\{G(X^{(i)}) > t\}}(X^{(i)})$\\
\tcc{For every subsets of inputs}
\For{$A \in \mathcal{P}_d$}{
    \tcc{Sample from the marginal distribution}
	$(X_A^{(1)}, \dots, X_A^{(N_v)}) \leftarrow \sM(A, N_v)$\\
	\tcc{For every element of the marginal distribution sample}
	\For{$i = 1, \dots, N_v$}{
	    \tcc{Sample from the conditional distribution given the element of the marginal distribution}
		$(\widetilde{X}_i^{(1)}, \dots, \widetilde{X}_i^{(N_p)})\leftarrow \sC(\overline{A}, N_p, X_A^{(i)})$
	}
	\tcc{Compute the conditional element}
	$\widehat{\textrm{T-S}}_A \leftarrow \frac{1}{N_v-1} \sum_{i=1}^{N_v} \left( \frac{1}{N_p} \sum_{j=1}^{N_p} \mathds{1}_{\mathcal{F}_t}(\widetilde{X}_i^{(j)},X_A^{(i)})  - \widehat{p}_t^Y  \right)^2 \times \frac{1}{\widehat{p}_t^Y (1-\widehat{p}_t^Y )}$
}

\tcc{Aggregation step}
\For{$j=1,\dots,d$}{
	$\widehat{\textrm{T-Sh}}_{j, \textrm{MC}} \leftarrow 0$\\
	\For{$A \subset \{-j\}$}{
        \tcc{Apply the Shapley weights to every computed increments}
		$\widehat{\textrm{T-Sh}}_{j, \textrm{MC}} += \frac{1}{d} {d-1 \choose |A|}^{-1} \left( \widehat{\textrm{T-S}}_{A \cup \{j\}} - \widehat{\textrm{T-S}}_{A} \right) $\\
	}
}
\caption{Target Shapley effects estimation by a Monte Carlo procedure.}
\label{alg:mc_est}
\end{algorithm}

Another algorithm has been proposed in \cite{song_shapley_2016}, by leveraging an equivalent definition of the Shapley \rev{allocations}, as an arithmetic mean over all the $d!$ permutations of $\{1, \dots,d \}$. 
In the same fashion as in \myeqref{eq:LMGperm}, it writes:
\begin{equation}
    \widehat{\textrm{T-Sh}}_{j, \textrm{MC}} = \frac{1}{m} \sum_{\substack{\pi \in \mbox{\scriptsize permutations}\\ \mbox{\scriptsize of } \{1,\ldots,d\}}} \left( \widehat{\textrm{T-S}}_{v \cup \{j\}, \textrm{MC}} - \widehat{\textrm{T-S}}_{v, \textrm{MC}} \right)
    \label{eq:shap_perm}
\end{equation}
with $v$ being the indices before $j$ in the order $\pi$. 
In \myeqref{eq:shap_perm}, the sum is not performed over all the permutations of $\{1,\ldots,d\}$ but only on $m$ randomly chosen permutations.
By sampling $m<d!$ permutations, one can drive the computational cost of this algorithm to $(N + m \times (d - 1) \times N_v \times N_p)$ calls to $G(\cdot)$, for a less precise, but still convergent estimator.

\subsection{Given-data estimation using a nearest-neighbor procedure}
\label{subsec:knn}
A ``given-data'' estimation method has been introduced by \cite{broto_variance_2020} to \rev{estimate} the Shapley effects. This method can be seen as an extension of the Monte Carlo estimator when only a single \iid input-output sample is available. This method is appropriate when the input distributions are not known or when the numerical model $G(\cdot)$ is \rev{not available anymore}. The main idea behind this method is to replace the exact samples from the conditional distributions $P_{X_{\overline{A}}|X_A}$ by approximated ones based on a non-parametric nearest-neighbor procedure.

Let $\left( X^{(1)}, \dots, X^{(N)}\right)$ be an \iid sample of the inputs $X$ and $A \in \mathcal{P}_d \setminus \{ \emptyset, \{1:d\}\}$. Let $k_N^A(l,n)$ be the index such that $X^{\left( k_N^A(l,n) \right)}_A$ is the $n$-th closest element to $X^{(l)}_A$ in $\left( X_A^{(1)}, \dots, X_A^{(N)}\right)$. Note that, if two observations are at an equal distance from $X^{(l)}_A$, then one of the two is uniformly randomly selected. Finally, one can define an estimator of the equivalent cost function defined in \myeqref{eq:alt_costfun}:
\begin{equation}
\widehat{\textrm{T-E}}_{A, \text{KNN}} =  \frac{\sum_{l=1}^N \left( \frac{1}{N_s -1} \sum_{i=1}^{N_s} \left[ \mathds{1}_{\mathcal{F}_t}\left(X^{\left( k^{\overline{A}}_N(l, i) \right)}\right)- \frac{1}{N_s} \sum_{h=1}^{N_s} \mathds{1}_{\mathcal{F}_t}\left(X^{\left( k^{\overline{A}}_N(l, h) \right)}\right)\right]^2\right)}{N\widehat{p}_t^Y(1-\widehat{p}_t^Y)}.
\label{eq:l2_knn}
\end{equation}

Under some mild assumptions, \cite{broto_variance_2020} showed that this estimator does asymptotically converge towards $\textrm{T-E}_A$. With estimates for the conditional elements, one can then define the following plug-in estimator:
\begin{equation}
\widehat{\textrm{T-Sh}}_{j, \text{KNN}} = \frac{1}{d} \sum_{A \subset \{-j\}} {d-1 \choose |A|}^{-1} \left( \widehat{\textrm{T-E}}_{A\cup\{j\}, \text{KNN}} - \widehat{\textrm{T-E}}_{A, \text{KNN}} \right)
\label{eq:l2tse_knn}
\end{equation}
where $\widehat{p}_t^Y$ is the empirical mean of $\mathds{1}_{\mathcal{F}_t}(X)$ on the \iid sample. \myalgref{alg:knn_est} represents the procedure for this given-data estimator. \rev{Its empirical convergence \wrt the sample size is illustrated in \ref{app:empconv_knn}}. This method is less computationally expensive (in terms of model evaluations) compared to the Monte Carlo sampling-based method, since no additional model evaluation, other than the ones in the \iid sample, is required in order to produce estimates of the target Shapley effects. Since the samples of the conditional and marginal distributions are approximated by a non-parametric procedure, this method also reduces the possible input modeling error (\eg in the context of ill-defined input distributions), at the cost of less accurate estimates. Another constraint is due to the fact that the input-output sample has to be \iid which prevents it from being used, for instance, in advanced \rev{orthogonal} designs of computer experiments.

\begin{algorithm}[!ht]
\SetKwFunction{knn}{KNN}

\KwIn{$X,Y,t$}
\KwOut{$(\widehat{\textrm{T-Sh}}_{j, \textrm{KNN}})_{j=1,\dots,d}$}
\tcc{Estimate the failure probability}
$\widehat{p}_t^Y \leftarrow \frac{1}{N} \sum_{i=1}^N \mathds{1}_{\{G(X^{(i)}) > t\}}(X^{(i)})$\\
\tcc{For every subsets of inputs}
\For{$A \in \mathcal{P}_d$}{
    \tcc{Sample of X\_A}
    $X_A \leftarrow (X_i^{(j)})_{j=1,\dots,n}^{i \in A}$\\
	\For{$i = 1, \dots, N$}{
	    \tcc{For each row i of X\_A, find the N\_s nearest rows in X}
		$(\widetilde{X}_{i}^{A, (j)})_{j=1, \dots, N_s}\leftarrow \knn(X_A^{(i)}, X, N_s)$
	}
	\tcc{Compute the conditional element}
	$\widehat{\textrm{T-E}}_A \leftarrow \sum_{l=1}^N \left( \frac{1}{N_s -1} \sum_{i=1}^{N_s} \left[ \mathds{1}_{\mathcal{F}_t}\left(\widetilde{X}_{l}^{A, (i)}\right)- \frac{1}{N_s} \sum_{h=1}^{N_s} \mathds{1}_{\mathcal{F}_t}\left(\widetilde{X}_{l}^{A, (h)}\right)\right]^2\right)\times \left(N\widehat{p}_t^Y(1-\widehat{p}_t^Y)\right)^{-1}$
}

\tcc{Aggregation step}
\For{$j=1,\dots,d$}{
	$\widehat{\textrm{T-Sh}}_{j, \textrm{MC}} \leftarrow 0$\\
	\For{$A \subset \{-j\}$}{
        \tcc{Apply the Shapley weights to every computed increments}
		$\widehat{\textrm{T-Sh}}_{j, \textrm{MC}} += \frac{1}{d} {d-1 \choose |A|}^{-1} \left( \widehat{\textrm{T-E}}_{A \cup \{j\}} - \widehat{\textrm{T-E}}_{A} \right) $\\
	}
}
\caption{Target Shapley effects estimation by a nearest-neighbor procedure.}
\label{alg:knn_est}
\end{algorithm}

In \cite{broto_variance_2020}, a random permutation algorithm, homologous to \myeqref{eq:shap_perm}, has been developed, which allows for reducing the overall complexity of the method, which, for the sake of conciseness, is not developed in this paper.

\subsection{Software and reproducibility of results}

The algorithms described in the preceding subsections have been implemented in the \texttt{sensitivity} R package \citep{ioodav20}. More precisely, the \texttt{shapleyPermEx()} (sampling-based algorithm) and \texttt{sobolshap\_knn()} (given-data algorithm) functions can be directly used for the estimation of the target Shapley effects. In the applications of Section \ref{sec:appli}, only the \texttt{sobolshap\_knn()} function is used for numerical tractability. \ref{sec:min_code} provides some minimal code examples for the implementation of the \rev{Monte Carlo (see, Section \ref{sec:mc_estim}) and nearest neighbors estimation procedure (see, Section \ref{subsec:knn})}, along with their random permutation variants.

All further results can be accessed on a \href{https://gitlab.com/milidris/review_l2tse}{GitLab}\footnote{https://gitlab.com/milidris/review\_l2tse} \rev{repository}, along with the data used in the following sections. R code files are available, with explicit code, along with all custom-made functions, in order to reproduce the \rev{results} presented in this paper. The procedures for the theoretical approximations of Section \ref{sec:toycases} are made available, along with the data-simulation functions for the flood case in Subsection \ref{subsec:flood_case}. The two datasets used for Subsection \ref{subsec:covid} are also available. Finally, all the figures can be reproduced by simply re-running the different \texttt{RMarkdown} files in the aforementioned GitLab repository.
\section{Analytical results using a linear model with Gaussian inputs }
\label{sec:toycases}

To illustrate the behavior of the proposed indices, a first toy-case involving a linear model and multivariate Gaussian inputs is presented. \rev{In this setting,} analytical results can be derived for the marginal distributions of all subsets of inputs, their conditional distribution, and the distribution of the output given a subset of inputs. \rev{Subsequently}, analytical formulas can be obtained for both the target Sobol' indices and the target Shapley effects.

Let $(\beta_0, \beta)=(\beta_0, \beta_1,\dots ,\beta_d)^{\top} \in \mathbb{R}^{d+1}$, $\mu=(\mu_1, \dots, \mu_d)^{\top} \in \mathbb{R}^d$ and $\Sigma \in \mathcal{M}_{d}(\mathbb{R})$ a full-rank symmetric $(d \times d)$ matrix. Assume that $X \sim \mathcal{N}_d\left( \mu, \Sigma \right)$, and that the model output writes
\begin{equation}
Y = \beta_0 + \beta^{\top} X.
\end{equation}
Then, one has $Y \sim \mathcal{N}\left( \beta_0 + \beta^{\top} \mu , \beta^{\top} \Sigma \beta\right)$ and, for any $A \in \mathcal{P}_d$, $(Y | X_A = x_A) \sim \mathcal{N} \left( \widetilde{\mu}_A , \widetilde{\Sigma}_A \right)$
with
$$\widetilde{\mu}_A = \beta_0 + \beta_A^{\top} x_a + \beta^{\top}_{\overline{A}}(\mu_{\overline{A}} + \Sigma_{A,12} \Sigma_{A,22}^{-1}(x_a - \mu_A)), \quad \widetilde{\Sigma}_A = \beta_{\overline{A}}^{\top} (\Sigma_{A,11} - \Sigma_{A,12}\Sigma_{A,22}^{-1}\Sigma_{A,21})\beta_{\overline{A}}.$$
Moreover, one also \rev{can} recall that
$$(X_{\overline{A}}, X_A)^{\top} \sim \mathcal{N}_d \left( \begin{pmatrix}\mu_{\overline{A}}\\ \mu_A \end{pmatrix}, \Sigma_A = \begin{pmatrix}\Sigma_{A,11}& \Sigma_{A,12}\\ \Sigma_{A,21} &\Sigma_{A,22}\end{pmatrix} \right) $$
with the partitions of $\Sigma_A$ having sizes $ \begin{pmatrix}(d-|A|) \times (d-|A|) & (d-|A|) \times |A| \\ |A| \times (d-|A|) & |A| \times |A| \end{pmatrix}$.
\rev{Evaluating these} results \rev{requires} some numerical approximations of the theoretical values of $\textrm{T-Sh}_j$ for all $j=1,\dots,d$. \rev{This has been achieved by} using standard multidimensional integration tools, \rev{and more specifically,} the function \texttt{adaptIntegrate()} from the \href{https://cran.r-project.org/web/packages/cubature/index.html}{\texttt{cubature}} package of the R software has been used, with a fixed error tolerance set to $10^{-8}$. This \rev{allowed the study of simple toy-cases} in order to validate the behavior of the target Shapley effects.

In the following, \rev{the inputs are first assumed to be independent, and are studied} \wrt the threshold $t$. Then, a \rev{toy-}case involving linear correlation between inputs (driven by \rev{a} parameter $\rho$) is studied. Finally, a last \rev{toy-}case \rev{aims at studying the proposed indices' behavior in the presence of} an exogenous input.

\subsection{Independent standard Gaussian inputs}

The first toy-case can be specified by:
\begin{equation}
\label{eq:LinGauIndep}
\begin{pmatrix}X_1 \\ X_2 \\ X_3 \end{pmatrix} \sim \mathcal{N}_3 \left( \begin{pmatrix}0\\0\\0 \end{pmatrix}, \begin{pmatrix}1 &0 &0 \\ 0&1&0 \\ 0&0&1 \end{pmatrix}\right), \quad Y = \sum_{i=1}^3 X_i.
\end{equation}
In this \del{scenario}\rev{case}, the three inputs are equally important in terms of defining $Y$, but they should also be equally important for the variable of interest $\mathds{1}_{\mathcal{F}_t}(X)$, as assessed by the target Sobol' indices defined as in \myeqref{eq:tsa_sobol}.

\begin{figure}[!ht]
\centering
\includegraphics[width=\textwidth]{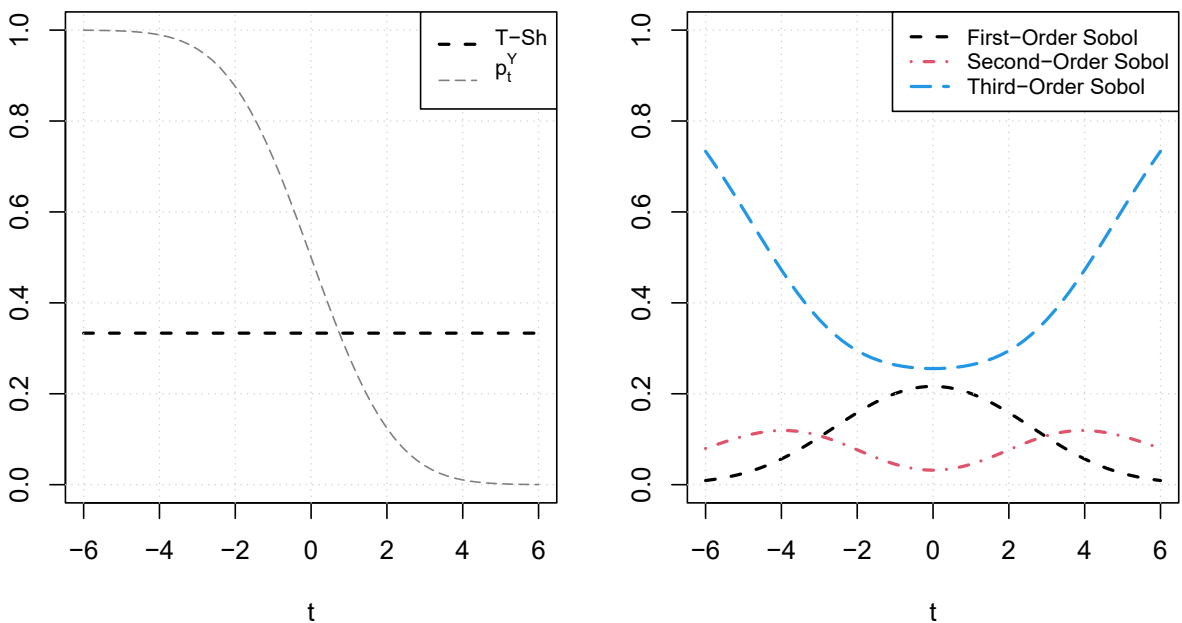}
\caption{Target Shapley effects (left) and target Sobol' indices (right) for the linear model with standard independent multivariate Gaussian inputs, \wrt $t$.}
\label{fig:LinGauIndepBal}
\end{figure}

From \citet{Li_Lu_SS_2012} and \citet{lemaitre_analyse_2014}, one can easily deduce that the first-order (FO) target \rev{closed} Sobol' indices are all equal to each other. Thus, one has:
\rev{\begin{equation}
\label{eq:firstorder_eq}
\textrm{T-S}_{\text{FO}} \overset{\text{def}}{=}\textrm{T-S}_1^{\ell^2}=\textrm{T-S}_2^{\ell^2}=\textrm{T-S}_3^{\ell^2} = \frac{\mathbb{V}\left( \Phi \left(\frac{t-X}{\sqrt{2}}\right)\right)}{\mathbb{V}\left( \mathds{1}_{\mathcal{F}_t}(X)\right)},
\end{equation}}
while the second-order (SO) target \rev{closed} Sobol' indices are given by:
\rev{\begin{equation}
\label{eq:secondorder_eq}
\textrm{T-S}_{\text{SO}} \overset{\text{def}}{=}\textrm{T-S}_{\{1,2\}}^{\ell^2}=\textrm{T-S}_{\{1,3\}}^{\ell^2}=\textrm{T-S}_{\{2,3\}}^{\ell^2}=\frac{\mathbb{V}\left( \Phi \left(t-X'\right)\right)}{\mathbb{V}\left( \mathds{1}_{\mathcal{F}_t}(X)\right)}
\end{equation}}
where $\Phi(.)$ is the standard Gaussian cumulative distribution function (cdf), $X\sim\mathcal{N}(0,1)$ and $X'\sim\mathcal{N}(0,2)$. Finally, one can also show that the third-order (TO) target \rev{closed} Sobol' indices are equal to:
\rev{\begin{equation}
\label{eq:thirdorder_eq}
\textrm{T-S}_{\text{TO}} \overset{\text{def}}{=} \textrm{T-S}_{\{1,2,3\}}^{\ell^2} = 1.
\end{equation}}
From Eqs. (\ref{eq:firstorder_eq}), (\ref{eq:secondorder_eq}), and (\ref{eq:thirdorder_eq}), and from Property \ref{prop:dtse_sum1}, one can deduce that:
\begin{equation}
\label{eq:l2tse_equal}
\textrm{T-Sh}_1 = \textrm{T-Sh}_2 = \textrm{T-Sh}_3 = \frac{1}{3} .
\end{equation}
Additionally, as the inputs are independent, interpreting the \rev{original target Sobol' indices (i.e., \myeqref{eq:original_tsasobol})} is meaningful, and they are equal to:
\rev{\begin{align}
\textrm{T-S}_i &= \textrm{T-S}_{\textrm{FO}}, \forall i \in \{1,2,3\} \\
\textrm{T-S}_{\{i,j\}} &= \textrm{T-S}_{\text{SO}} - 2 \textrm{T-S}_{\text{FO}}, \forall i,j \in \{1,2,3\}, i \not = j \\
\textrm{T-S}_{\{1,2,3\}} &= \textrm{T-S}_{\text{SO}} -  3 \left( \textrm{T-S}_{\text{FO}} + \textrm{T-S}_{\text{SO}}\right).
\end{align}}
\rev{The target Sobol' indices are illustrated in \myfigref{fig:LinGauIndepBal} (right)}. One can remark that, \rev{focusing on the indicator variable of interest} \del{studying the indicator variable} $\mathds{1}_{\mathcal{F}_t}(X)$ instead of the model output $Y$ leads to interaction effects between the inputs, \rev{as outlined in Section \ref{subsec:tsa}. The target Shapley effects, however,} remain constant for all threshold values $t$. Such a behavior \rev{is} expected\rev{:} it highlights the fact the target Shapley effects do not report the interaction effects as \rev{the} target Sobol' indices would. \rev{The proposed indices rather summarize (in the sense of the Shapley values allocation) the target Sobol' indices into a single index. Their goal is not to report on the ``types of effects'' (\ie correlation or interaction), but rather provide a global index which sums up each input's importance.}

\subsection{Correlated Gaussian inputs with unit variance}

The behavior of the target Shapley effects are now studied when a linear dependence is added to the inputs. Since Property \ref{prop:dtse_sum1} still holds without any condition on the dependence structure on the input variables, these indices remain \rev{interpretable as a percentage of the output's variance}. 
The following model is \rev{studied}:
\begin{equation}
\label{eq:LinGauDep}
\begin{pmatrix}X_1 \\ X_2 \\ X_3 \end{pmatrix} \sim \mathcal{N}_3 \left( \begin{pmatrix}0\\0\\0 \end{pmatrix}, \begin{pmatrix}1 &0 &0 \\ 0&1&\rho \\ 0&\rho&1 \end{pmatrix}\right), \quad Y = \sum_{i=1}^3 X_i.
\end{equation}
where $-1<\rho<1$. In this scenario one has:
\rev{\begin{align}
\label{eq:fo_dep}
\textrm{T-S}_1^{\ell^2} &= \frac{\mathbb{V}\left(  \Phi \left( \frac{t-X}{\sqrt{2(1+\rho)}} \right) \right)}{\mathbb{V}\left(\mathds{1}_{\mathcal{F}_t}(X) \right)},\\
\textrm{T-S}_2^{\ell^2}&=\textrm{T-S}^{\ell^2}_3 = \frac{\mathbb{V}\left(  \Phi \left( \frac{t-X(1+\rho)}{\sqrt{2-\rho^2}} \right) \right)}{\mathbb{V}\left(\mathds{1}_{\mathcal{F}_t}(X) \right)},\\
\textrm{T-S}_{\{1,2\}}^{\ell^2}&=\textrm{T-S}^{\ell^2}_{\{1,3\}}=\frac{\mathbb{V}\left(  \Phi \left( \frac{t-X'}{\sqrt{1-\rho^2}}\right) \right)}{\mathbb{V}\left(\mathds{1}_{\mathcal{F}_t}(X) \right)},\\
\textrm{T-S}_{\{2,3\}}^{\ell^2}&=\frac{\mathbb{V}\left(  \Phi \left( t-X''\right) \right)}{\mathbb{V}\left(\mathds{1}_{\mathcal{F}_t}(X) \right)}
\end{align}}
where $X\sim \mathcal{N}(0, 1)$, $X'\sim \mathcal{N}(0, 1+(1+\rho)^2)$ and $X''\sim \mathcal{N}\left(0, 2(1+\rho)\right)$.

\begin{figure}[!ht]
\includegraphics[width=\textwidth]{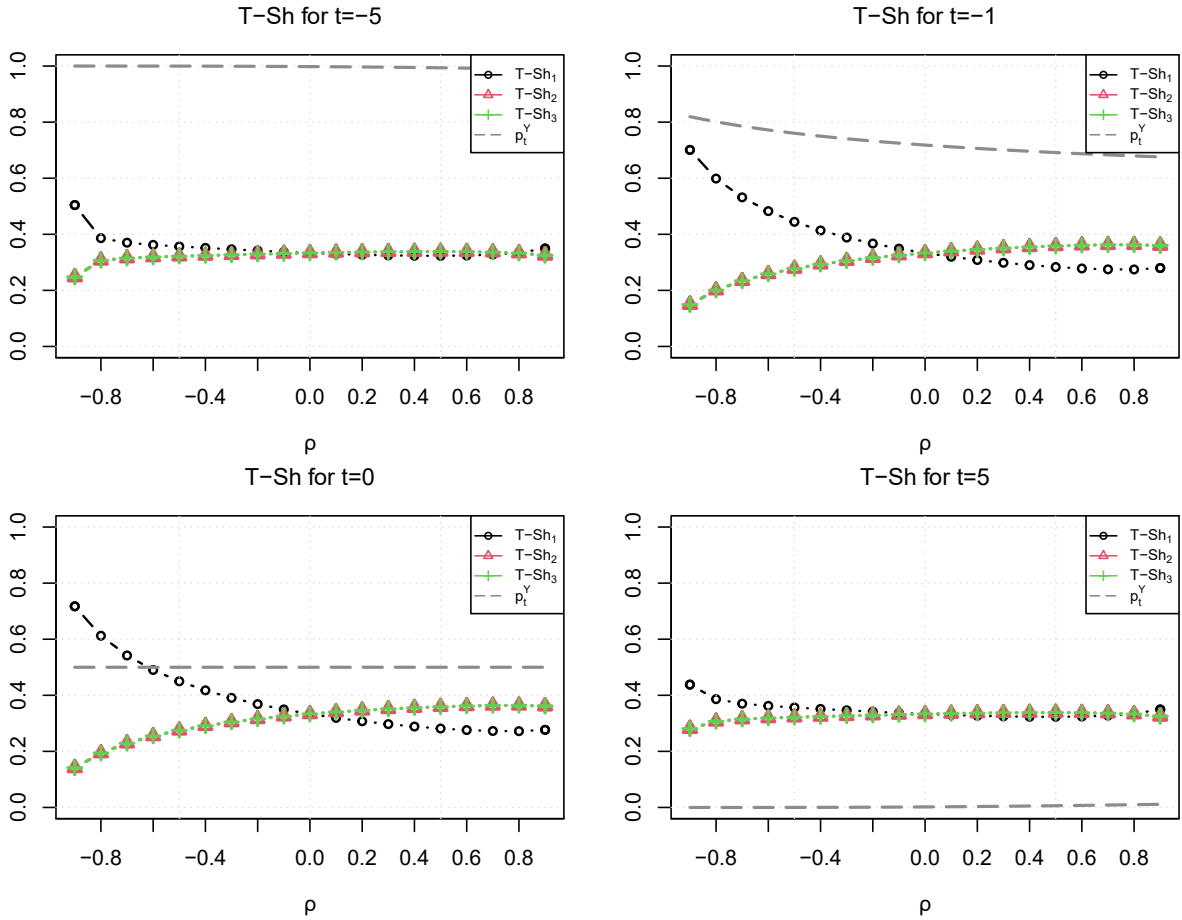}
\caption{Evolution, \wrt $\rho$ and for various threshold values, of the target Shapley effects of correlated Gaussian standard inputs in a linear model.}
\label{fig:l2tse_LinGauDepRho}
\end{figure}

\begin{figure}[!ht]
\centering
\includegraphics[width=0.8\textwidth]{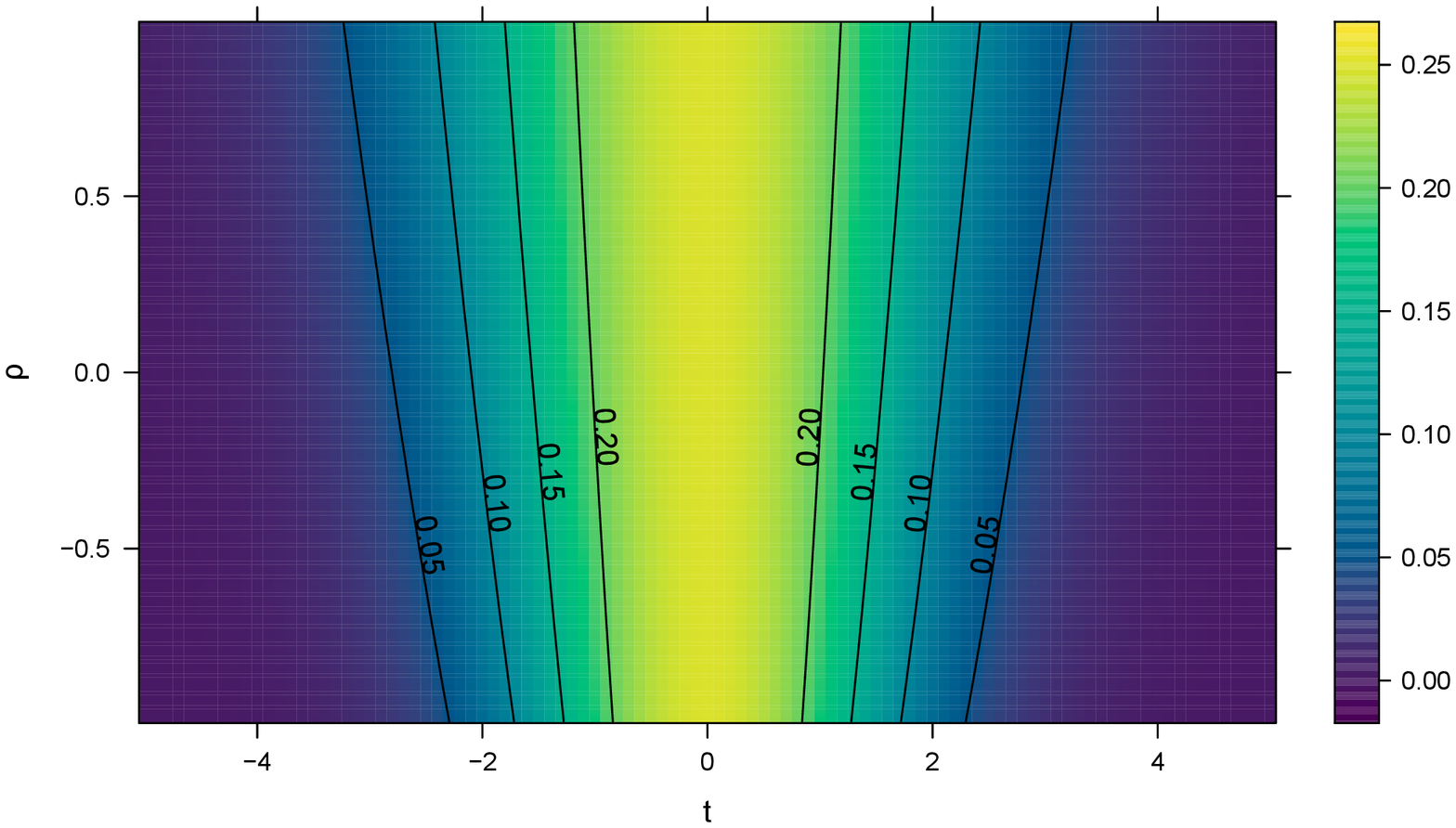}
\caption{Variance of $\mathds{1}_{\mathcal{F}_t}(X)$ \wrt $\rho$ and $t$ for correlated Gaussian standard inputs with a linear model.}
\label{fig:var1f_GauLinDep}
\end{figure}

From these results, one can directly remark that $\textrm{T-Sh}_2 = \textrm{T-Sh}_3$.
Note that the values of the target Shapley effects can also be obtained by combinations of target Sobol' indices (see \myeqref{eq:l2tse}). These results are illustrated in \myfigref{fig:l2tse_LinGauDepRho}.
For fixed threshold values $t$, the target Shapley effects of the correlated inputs $X_2$ and $X_3$ increases when $\rho$ increases. This is an expected behavior since, in this case:
\begin{equation}
\label{eq:varval_LinGauDep}
\mathbb{V}\left( \mathds{1}_{\mathcal{F}_t}(X)\right) = \Phi\left( \frac{t}{\sqrt{3+2\rho}}\right) \left(1- \Phi\left( \frac{t}{\sqrt{3+2\rho}}\right) \right),
\end{equation}
and subsequently, for a fixed $t$, the variance of the variable of interest will grow with $\rho$, as illustrated in \myfigref{fig:var1f_GauLinDep}. This increase in variance due to the correlation between $X_2$ and $X_3$ is then \del{allocated}\rev{attributed} through $\textrm{T-Sh}_2$ and $\textrm{T-Sh}_3$, which increase with $\rho$. On the other hand, $\textrm{T-Sh}_1$ decreases accordingly, to accommodate Property \ref{prop:dtse_sum1}.

In \myfigref{fig:l2tse_LinGauDepRho}, the behavior of the indices \wrt $\rho$ \rev{is illustrated.} $\textrm{T-Sh}_1$ is predominantly above $\textrm{T-Sh}_2$ and $\textrm{T-Sh}_3$ when $\rho$ is negative, and below when it is positive. This can be explained by the fact that, $X_2$ and $X_3$ cancel each other out when their correlation is negative, thus lowering the value of $\textrm{T-S}_{\{2,3\}}$ below $\textrm{T-S}_{\{1,2\}}$ and $\textrm{T-S}_{\{1,3\}}$, automatically increasing $\textrm{T-Sh}_1$ in accordance to Property \ref{prop:dtse_sum1}. \rev{On} the other hand, for positive values of $\rho$, $\textrm{T-S}_{\{2,3\}}$ is higher than $\textrm{T-S}_{\{1,2\}}$ and $\textrm{T-S}_{\{1,3\}}$, which in turn corresponds to $\textrm{T-Sh}_1$ being lower than $\textrm{T-Sh}_2=\textrm{T-Sh}_3$.

\subsection{Quantifying the importance of an exogenous input in the Gaussian setting}
\label{sec:exovar}
\begin{figure}[!ht]
\centering
\includegraphics[width=0.8\textwidth]{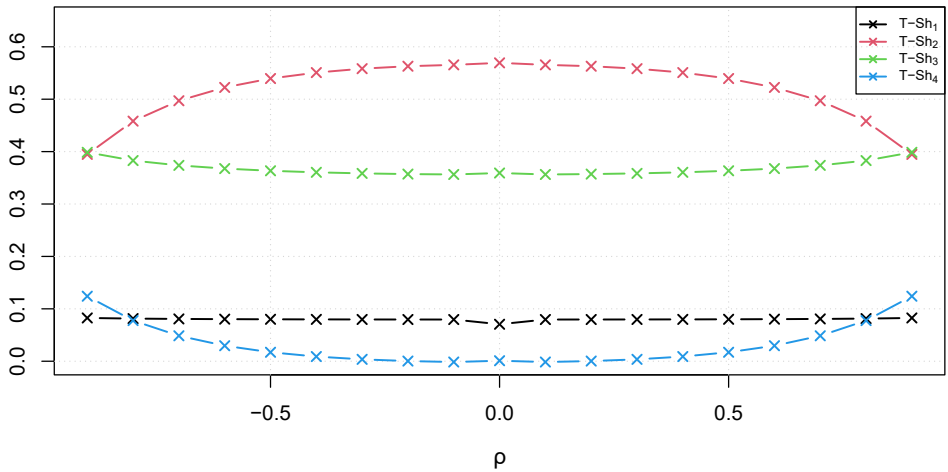}
\caption{Target Shapley effects for the Gaussian Linear model with an exogenous input, \wrt $\rho$.}
\label{fig:exo_variable_fig}
\end{figure}

In this \del{use }\rev{toy-}case, inspired by \cite{lemaitre_analyse_2014}, the following model is considered:
\begin{equation}
\begin{pmatrix}X_1\\X_2\\X_3\\X_4 \end{pmatrix} \sim \mathcal{N}_4  \left( \begin{pmatrix}0\\0\\0\\0 \end{pmatrix}, \begin{pmatrix}1 &0 &0&0 \\ 0&1&0&\rho \\ 0&0&1&0\\0&\rho&0&1 \end{pmatrix}\right), \quad Y =X_1 + 6X_2 + 4X_3
\label{eq:model_exo}
\end{equation}
where $X_4$ is an exogenous input, but correlated to $X_2$, which \rev{is} the most important variable in terms of variance \rev{contribution, due to its higher linear coefficient}. The threshold is fixed at $t=16$. This scenario \rev{allows} for the verification of \rev{how} the target Shapley effects \rev{attribute the importance of} $X_4$ \rev{which is correlated} with an endogenous input, even though it does not appear in the model. \rev{In the} results given in \myfigref{fig:exo_variable_fig}, one can remark that $\textrm{T-Sh}_4$ increases when $\rho$ \rev{approaches} either $1$ or $-1$, despite the fact that it \rev{has no direct causal effect on} the model $G(\cdot)$.

\rev{
\subsection{Discussion on causal relationship assessment}
As outlined in Subsection \ref{sec:exovar}, one can remark that, in the context of highly correlated inputs, the proposed indices fail to provide insights on the causal relationships of an input on the model's output: an exogenous input may receive a share of the output's variance.
This behavior is intrinsically due to the Shapley values allocation method, and has been highlighted in \cite{iooss_shapley_2019} in the case of the Shapley effects. This particularity makes the interpretation of the Shapley effects, and subsequently the target Shapley effects, quite delicate. A prior investigation of the correlation structure of the data, through, for instance, estimated correlation matrices or any tool dedicated to multicollinearity diagnostics such as the variance inflation factor \citep{foxmon92} is strongly advised, along with an input validation process, ensuring the absence of exogenous inputs in the TSA study. Other allocations systems, such as the \textit{proportional values} \citep{ortmann_proportional_2000}, are purposefully designed in order to highlight causal effects \citep{feldman_relative_2005}, but fall out of the scope of this document and are not developed further.}
\section{Applications}
\label{sec:appli}
 
In this section, two models related to real\rev{-world} phenomena which include dependent random inputs are studied in the context of TSA.
 
\subsection{A simplified flood model}
\label{subsec:flood_case}

The target Shapley effects are firstly computed on a simplified model of a river flood \citep{lemaitre_analyse_2014,ioolem15}.
\rev{This model's goal} is to simulate the behavior of a river's water level, \rev{and to compare it to a fixed dyke height}. After a strong simplification of the one-dimensional Saint-Venant equation (with uniform and constant flow rate), the maximal annual water level $h$ is modeled as:
\begin{equation}
h = \left( \frac{Q}{BK_s \sqrt{\frac{Z_m - Z_v}{L}}}\right)^{\frac{3}{5}},
\end{equation}
while the model output writes:
\begin{equation}
Y=Z_v + h.
\end{equation}
The \rev{six} inputs' \rev{probabilistic structure is} described in Table~\ref{tab:tabflood}. The problem is of dimension $d=6$. \rev{Under the TSA paradigm,} the variable of interest is $\mathds{1}_{\{ G(X) > t\}}(X)$ with $t$ \rev{representing} the dyke\rev{'s} height, fixed to $t=54.5~\textrm{m}$. The reference failure probability \rev{(see, \myeqref{eq:tsa_qoi})}, computed here with a Monte Carlo sample of large size (here $10^7$ samples) is equal to $p_t^Y = 4.5 \times 10^{-3}$.

\begin{table}[!ht]
\centering
\begin{tabular}{|l|l|l|l|}
\hline
Input & Description & Unit & Distribution\\
\hline
\hline
$Q$ & maximal annual flow rate & $\textrm{m}^3.\textrm{s}^{-1}$ & Gumbel$(1013,558)$ truncated to $[500,3000]$\\
\hline
$K_s$ & Strickler friction coefficient &  $\textrm{-}$ & Normal$(30,7)$ truncated to $[15, +\infty)$\\
\hline
$Z_v$ & river downstream level  & $\textrm{m}$ & Triangular$(49,50,51)$\\
\hline
$Z_m$ & river upstream level & $\textrm{m}$ & Triangular$(54,55,56)$\\
\hline
$L$ & length of the river stretch & $\textrm{m}$ & Triangular$(4990,5000,5010)$\\
\hline
$B$ & river width & $\textrm{m}$ & Triangular$(295,300,305)$\\
\hline
$t$ & dyke height (threshold) &  $\textrm{m}$ & Fixed to $54.5$\\
\hline
\end{tabular}
\caption{Input variables and distributions for the flood model.}
\label{tab:tabflood}
\end{table}

\rev{In the same fashion as} in \cite{chastaing_generalized_2012}, three pairs of inputs are assumed to be linearly dependent: $Q$ and $K_s$ with $\rho(Q,K_s) = 0.5$, $Z_v$ and $Z_m$ with $\rho(Z_v,Z_m) = 0.3$, $L$ and $B$ with $\rho(L,B) = 0.3$. The aim of this use-case is to assess the \rev{relevance} of the target Shapley effects in a complex environment. \rev{In} \cite{chastaing_generalized_2012}, it \rev{is} shown that, from a GSA standpoint (using a generalized variance decomposition for dependent variables), the two most influential inputs on the annual water level are $Q$, the maximal annual flow rate, and $Z_v$, the river downstream level.

\begin{figure}[!ht]
\centering
\includegraphics[width=0.8\textwidth]{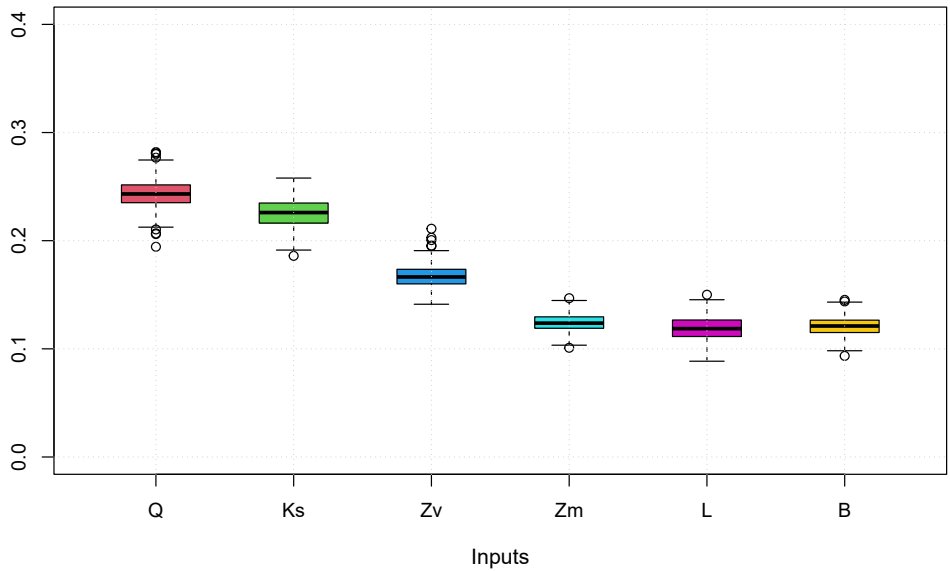}
\caption{Estimated target Shapley effects for the flood case.}
\label{fig:floodcase_est}
\end{figure}

\rev{An \iid} sample of $N=2\times 10^5$ input realizations is drawn (note that the \rev{linear} correlations are injected following the algorithm proposed by \cite{schumann_algo_2009}) \rev{which} leads to $N$ model \rev{evaluations}.
\myfigref{fig:floodcase_est} \rev{presents} the estimated target Shapley effects on \rev{this \iid sample}, using the nearest-neighbor procedure depicted in Subsection~\ref{subsec:knn} with an arbitrary number of neighbors set at $N_s=2$. $300$ repetitions of the simulation and the estimation procedure \rev{allow for the assessment of the estimation procedure's variance} (represented by boxplots in \myfigref{fig:floodcase_est}). \rev{One} can notice that $Q$ is granted an influence of $24.3\%$ ($\pm 1.3\%$), $K_s$ has $22.6\%$ ($\pm$1.3\%) and $Z_v$ around $16.7\%$ ($\pm 1\%$). The other inputs are attributed a share of around $12\%$. 
Compared to results obtained by \rev{from a} GSA \rev{standpoint,} without correlations \citep{ioolem15} and with correlations \citep{chastaing_generalized_2012}, these TSA results \rev{allow for granting} a much larger \rev{share} to $K_s$ and non-negligible effects to $Z_m$, $L$ and $B$. This was expected due to the interactions induced by the considered TSA variable of interest. This example illustrates the ability of the target Shapley effects to quantify the importance of input variables in a \rev{use-case} \rev{involving input correlation.}

\subsection{A COVID-19 epidemiological model}
\label{subsec:covid}

In 2020, the COVID-19 pandemic has raised important questions on the usefulness of epidemic modeling, especially on their ability to produce relevant insights to public policy decision makers. \cite{salbam20} have taken this example to insist on the essential use of GSA on such models, which claim to predict the potential consequences of intervention policies. A first study has been proposed by \cite{lubor20}, in the context of COVID-19 in Italy, to assess the sensitivity of important epidemiological model outcomes, such as the number of people being either quarantined, recovering, or dead due to COVID-19. Another GSA has been performed in \cite{veiga_basics_2020} in the French context of the first COVID-19 outbreak. By using data coming from this last analysis (thanks to the authors' agreement), the goal of this section is to demonstrate how TSA can help to characterize the influence of various uncertain parameters on a real-scale model.

\subsubsection{Description of the model and its inputs}
\begin{figure}[!ht]
\centering
\includegraphics[width=0.8\textwidth]{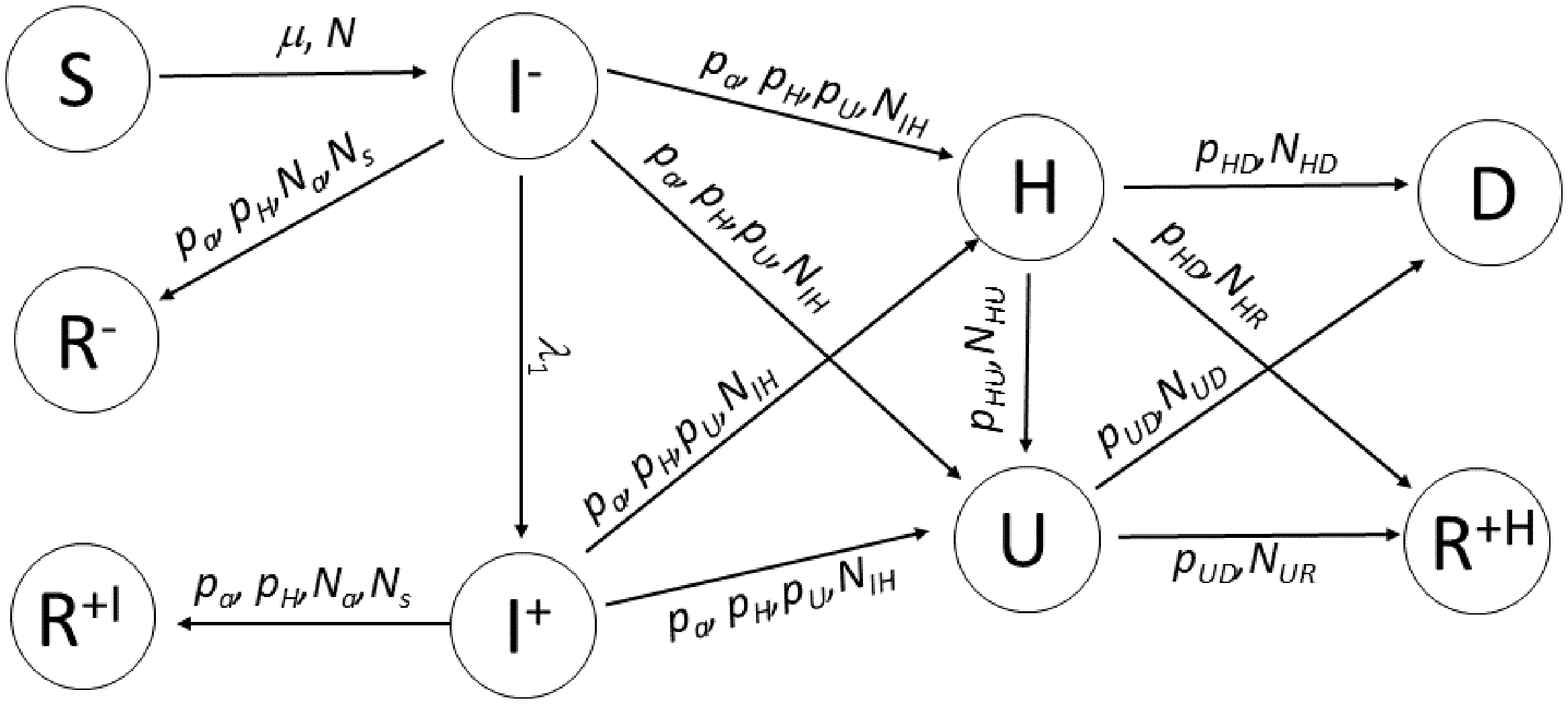}
\caption{Overview of the COVID-19 compartmental model. The different compartments contain the following individuals: susceptible (S), undetected infected (I$^-$), unhospitalized detected infected (I$^+$), recovered undetected (R$^-$), recovered unhospitalized detected (R$^{+I}$), hospitalized in severe illness (H), hospitalized in intensive care (U), recovered hospitalized (R$^{+H}$) and dead (D).}
\label{fig:covid19}
\end{figure}

The deterministic compartmental model developed in \cite{veiga_basics_2020} (also presented in \cite{dav20}) is representative of the COVID-19 French epidemic (from March to May) by taking into account the asymptomatic individuals, the testing strategies, the hospitalized individuals, and people admitted to Intensive Care Unit (ICU). Using several assumptions, it is based on a system of $10$ ordinary differential equations that can be fully retrieved in references \cite{dav20} and \cite{veiga_basics_2020}.
\rev{Each equation models path of individuals between different compartments (corresponding to their infectious and illness states), as shown in \myfigref{fig:covid19}.
These equations involve many input parameters and model the dynamic between the different compartments.}
Table~\ref{tab:prior_inputs} presents the $20$ \rev{continuous} input parameters with their prior distribution (chosen from literature studies), which form the inputs $X$, assumed to be independent between each other.

\begin{table}[!ht]
    \centering
    \begin{tabularx}{\textwidth}{|l|X|l|}
    \hline
         Input & Description & Prior distribution \\
         \hline
         \hline
         $p_a$ & Conditioned on being infected, the probability of having mild symptoms or no symptoms & $\mathcal{U}(0.5,0.9)$\\ 
         \hline
         $p_H$ &Conditioned on being mild/severely ill, the probability of needing hospitalization ($H$ or $U$) &$\mathcal{U}(0.15,0.2)$ \\ 
         \hline
         $p_U$ &Conditioned on hospitalisation, the probability of being admitted to ICU &$\mathcal{U}(0.15,0.2)$ \\ 
         \hline
         $p_{HD}$ &Conditioned on being hospitalized but not in ICU, the probability of dying &$\mathcal{U}(0.15,0.25)$ \\
         \hline
         $p_{UD}$ &Conditioned on being admitted to ICU, the probability of dying &$\mathcal{U}(0.2,0.3)$ \\ 
         \hline
         $N_a$ &If asymptomatic, number of days until recovery &$\mathcal{U}(8,12)$ \\ 
         \hline
         $N_s$ &If symptomatic, number of days until recovery without hospitalisation &$\mathcal{U}(8,12)$ \\ 
         \hline
         $N_{IH}$ &If severe symptomatic, number of days until hospitalization & $\mathcal{U}(8,12)$\\ 
         \hline
         $N_{HD}$ &If in $H$, number of days until death & $\mathcal{U}(15,20)$\\ 
         \hline
         $N_{UD}$ &If in ICU, number of days until death & $\mathcal{U}(8,12)$\\ 
         \hline
         $N_{HR}$ &If hospitalized but not in ICU, the number of days until recovery & $\mathcal{U}(15,25)$\\ 
         \hline
         $N_{UR}$ &If in ICU, number of days until recovery & $\mathcal{U}(15,25)$\\ 
         \hline
         $R_0$ &Basic reproduction number & $\mathcal{U}(3,3.5)$\\ 
         \hline
         $t_0$ &Starting date of epidemics (in 2020) & $\mathcal{U}(01/25,02/24)$\\ 
         \hline
         $\mu$ &Decaying rate for transmission (after social distanciation and lockdown) & $\mathcal{U}(0.03,0.08)$\\ 
         \hline
         $N$ &Date of effect of social distanciation and lockdown & $\mathcal{U}(20,50)$\\ 
         \hline
         $\lambda_1$ &Type-1 testing rate & $\mathcal{U}(1e-4,1e-3)$\\ 
         \hline
         $p_{HU}$ &Conditioned on being hospitalized in $H$, the probability of being admitted to ICU & $\mathcal{U}(0.15,0.2)$\\ 
         \hline
         $N_{HU}$ & If in $H$, number of days until ICU& $\mathcal{U}(1,10)$\\ 
         \hline
         $I^{-}_0$ &Number of infected undetected at the start of epidemics & $\mathcal{U}(1,100)$\\ 
         \hline
    \end{tabularx}
    \caption{Model inputs and their prior distribution. $H$ is the number of hospitalized individuals with severe symptoms. $U$ is the number of hospitalized individuals in ICU.}
    \label{tab:prior_inputs}
\end{table}

For the present study, our variable of interest, which is a particular model output, then writes
\begin{equation}
    \label{eq:model_prior}
    U_{\textrm{max}}^{\textrm{p}} = \underset{v \in \textrm{time range}}{\textrm{max }} \Bigl\{ U_v(X) \Bigr\}
\end{equation}
where $U_v$ is the the number of hospitalized patients in ICU at time $v$. 
Note that the ``$\textrm{p}$'' in $U_{\textrm{max}}^{\textrm{p}}$ stands for ``prior'' as this quantity corresponds to the variable of interest before any calibration \wrt the available data.

In \cite{veiga_basics_2020}, after a first screening step which allows for suppressing non-influential inputs, the model is calibrated on real data by using a Bayesian calibration technique.
After the analysis of this step, the selected remaining inputs are
\begin{equation}
\label{eq:inputs_post}
    X_{\textrm{sel}}=(p_a,N_a,N_s,R_0,t_0,\mu,N,I^{-}_0)^{\top}
\end{equation}
and their distributions are obtained from a sample given by the calibration process.
The non-influential inputs are fixed to their nominal values and the posterior variable of interest becomes
\begin{equation}
\label{eq:model_post}
U_{\textrm{max}}= \underset{v \in \textrm{time range}}{\textrm{max }} \Bigl\{ U_v(X_{\textrm{sel}}) \Bigr\}
\end{equation}
with $U_{\textrm{max}}$ being the maximum number of hospitalized people in ICU who need special medical care on the studied temporal range, and $U_v$ is the number of hospitalized patients in ICU at time $v$.

\subsubsection{Input importance for ICU bed shortage}

The central question of this study would be to determine which inputs influence the event of a country experiencing a shortage of ICU bed capacity during the time period. For that purpose, one can introduce a threshold \rev{$k$}, which represents the total number of ICU beds in the country, which is assumed to be constant during the studied time period. The new variable of interest would then be \rev{$\mathds{1}_{\left\{ U_{\textrm{max}}^{\textrm{P}} > k\right\}}(X)$} for the full compartmental model (preliminary study) and \rev{$\mathds{1}_{\{U_{\textrm{max}>k}\}}(X_{\textrm{sel}})$} for the model with selected inputs (post-calibration study). Two input-output samples of size $n=5000$ are available.
The first one (preliminary study) includes all the inputs following their prior distribution (see Table~\ref{tab:prior_inputs}) and the corresponding output $U^P_{\textrm{max}}$ of the compartmental model.
The second one (post-calibration study) is composed of a sample of $X_{\textrm{sel}}$ after the Bayesian calibration, and the corresponding output $U_{\textrm{max}}$ of the compartmental model with the non-selected inputs fixed to their nominal values.

\begin{figure}[!ht]
    \centering
    \includegraphics[width=\textwidth]{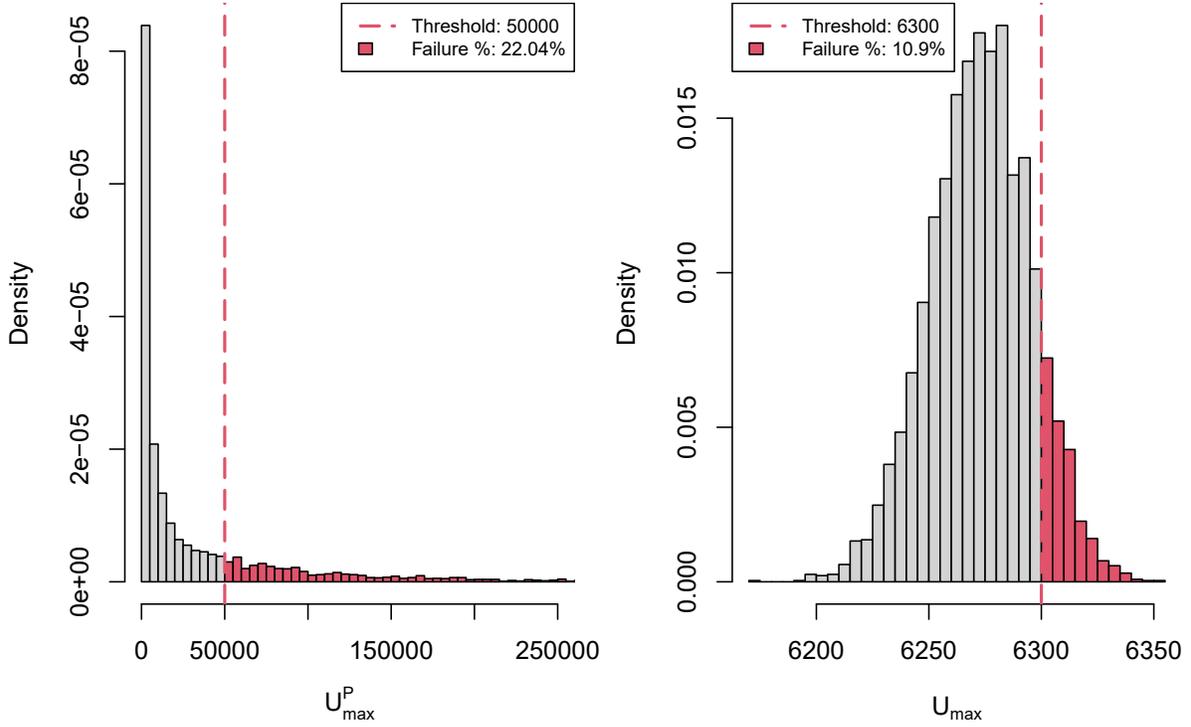}
    \caption{Illustration of thresholds on the histograms of $U^P_{\textrm{max}}$ (left) and $U_{\textrm{max}}$ (right).}
    \label{fig:hist_covid}
\end{figure}

Five different thresholds are studied on $U^P_{\textrm{max}}$: $5\cdot10^3$, $10^4$, $5\cdot10^4$, $10^5$ and $2\cdot10^5$, with respectively $58.1\%$, \rev{$47.7\%$}, $22\%$, $10.1\%$ and $2.2\%$ of the total output samples being in a failure state. This \del{allows to }illustrates the behavior of the target Shapley effects when the failure probability decreases. The threshold of $6300$ has been chosen for $U_{\textrm{max}}$, with $10.9\%$ of the total output samples being above this threshold. \myfigref{fig:hist_covid} illustrates two different thresholds, and the corresponding estimated failure probability on the histogram of both outputs.

\begin{figure}[!ht]
    \centering
    \includegraphics[width=0.8\textwidth]{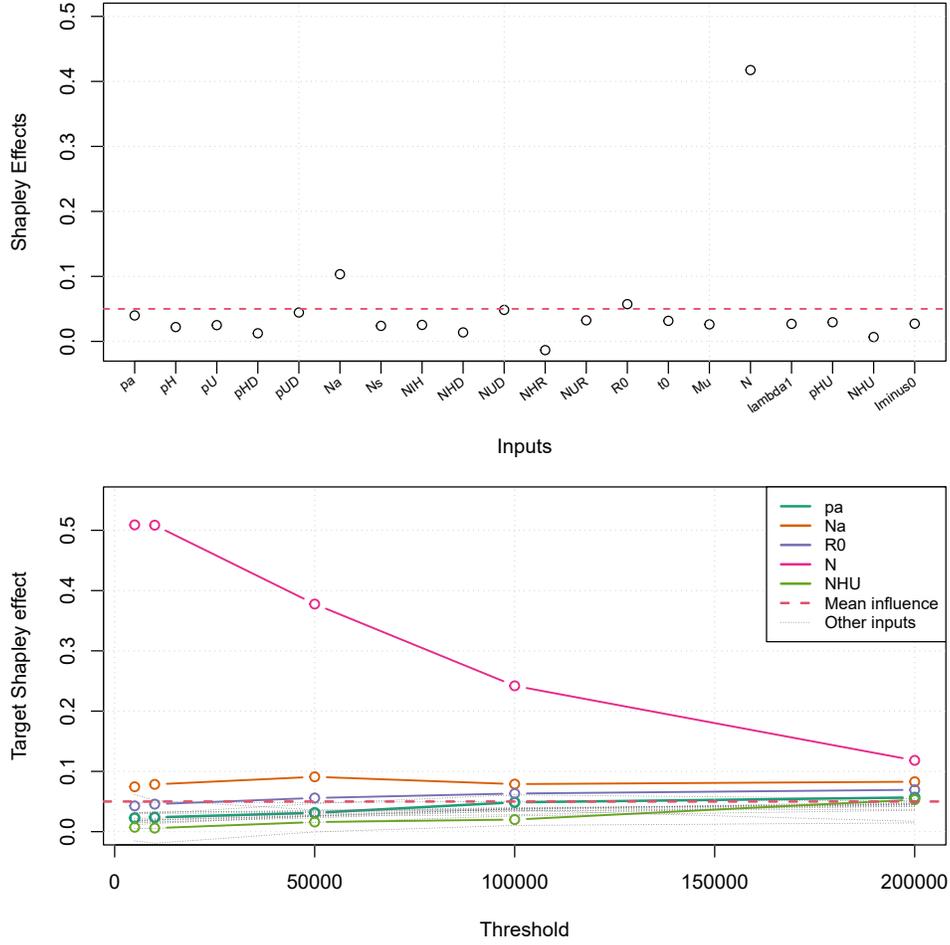}
    \caption{Shapley effects (top) and target Shapley effects for different thresholds (bottom) for $ U_{\textrm{max}}^P$.}
    \label{fig:l2tse_prior}
\end{figure}

The target Shapley effects have been estimated using a variant of the estimation scheme presented in Subsection~\ref{subsec:knn}, with a fixed number of random permutations of $10^3$, and with a number of neighbors set to $3$, following the rule of thumb guideline of \cite{broto_variance_2020}, due to the sheer complexity of this estimation algorithm. Since the compartmental model is deterministic, the target Shapley effects have been forced to sum up to one. \myfigref{fig:l2tse_prior} presents the main results for $U_{\textrm{max}}^P$, with the red dotted line being the average influence of an input, in the case of similar importance (\ie $\frac{1}{20}$). One can remark that for less restrictive thresholds (\ie threshold for which the failure probability is high), the input $N$, the effective date of lockdown/social distanciation measures, seem to be the most influential, reaching more than $50\%$ of the TSA variable of interest's variance. However, as soon as the threshold becomes more and more restrictive (\ie the failure probability decomes lower and lower), the effect of $N$ decreases, and the effects of the other inputs increase accordingly, in order to reach what seem to be an equilibrium at the value $\frac{1}{20}$. This behavior can be explained by two main reasons:
\begin{itemize}
    \item As outlined in Subsection~\ref{subsec:tsa}, the nature of a restrictive TSA variable of interest induces high interaction between the inputs;
    \item The Shapley allocation system, when applied to variance as a production value, redistributes the interaction effects equally between all inputs (there is no correlation between inputs in this prior study).
\end{itemize}
One can argue that, as soon as \rev{$k$} becomes very restrictive, the combined interaction effects outweighs the effect of $N$ itself, and since these effects are equally distributed among all the inputs, their share will tend to go towards $\frac{1}{20}$.

\begin{figure}[!ht]
    \centering
    \includegraphics[width=\textwidth]{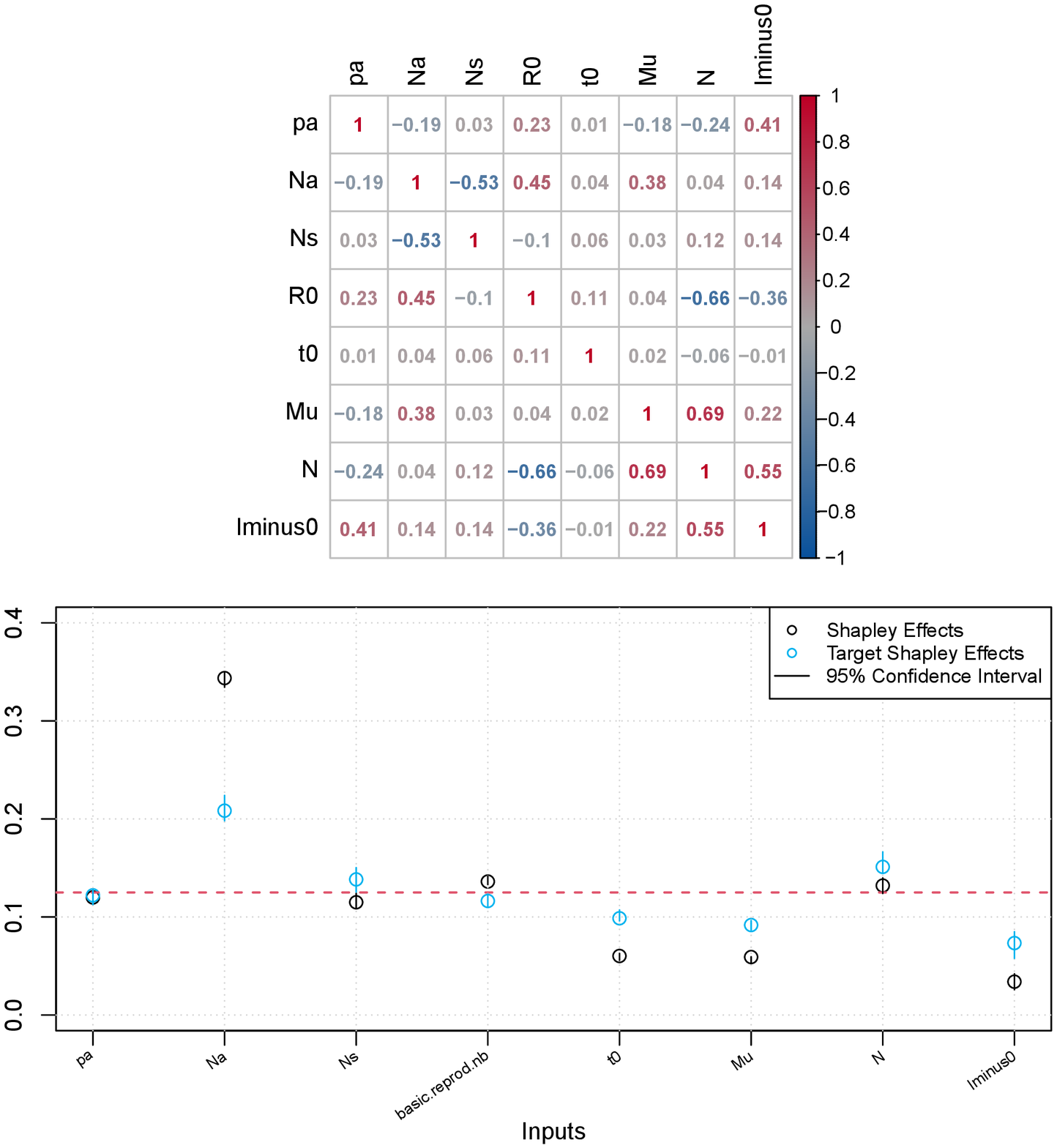}
    \caption{Input correlation matrix (top), Shapley effects for $U_{\textrm{max}}$ and target Shapley effects (bottom) for $\mathds{1}_{\{U_{\textrm{max}>t}\}}(X_{\textrm{sel}})$. The 95\% confidence intervals have been computed by uniformly selecting $80\%$ of the observations, for $100$ repetitions, without replacement.}
    \label{fig:corTshap_post}
\end{figure}

For the post-calibration study, some selected inputs $X_{\textrm{sel}}$ are linearly correlated (see \myfigref{fig:corTshap_post} - top\del{ left}). This is typically the case for $N$ and $\mu$, with an estimated correlation coefficient $\hat{\rho}(N,\mu)=0.69$, and for $R_0$ and $N$ with an estimated correlation coefficient of $\hat{\rho}(N,R_0)=-0.66$. This correlation structure does not allow for interpretable Sobol' indices, as outlined in Section~\ref{sec:gsa_dep}, which encourages the use of Shapley-inspired indices. The Shapley effects and the target Shapley effects of $X_{\textrm{sel}}$ for $U_{\textrm{max}}$ have been computed using the nearest-neighbor procedure, with a fixed number of neighbors of $3$, and forced to sum to one because of the deterministic nature of the model.

In \myfigref{fig:corTshap_post} (bottom), one can remark that $N_a$, the number of days until recovery, seem to be the most important input in explaining the number maximum number of ICU patients on the studied time range, with a Shapley effect of around $35\%$ of the output variance. The inputs $p_a$, $N_s$, $R_0$ and $N$ seem to present average effects, that is around $\frac{1}{8}$, while $t_0$, $\mu$ and $I^{-}_0$ seem to be less influential, with around $5\%$ of explained variance each.

However, focusing on the occurrence of a ICU bed shortage, one can remark that the target Shapley effect of $N_a$ is lower (around $22\%$), with the influence of $N$ being higher (around $15\%$) than their Shapley effects. Moreover, $t_0$, $\mu$ and $I^{-}_0$ present higher TSA effects, \ie slightly under $10\%$, due to the interaction induced by the indicator function. One can also remark that the influence of $N_s$ is higher than that of $R_0$ in the TSA setting, which was the inverse for the Shapley effects. This would indicate that $N_s$, the number of days until recovery for a symptomatic patient without hospitalization, has more influence on the event of a bed shortage than the basic reproducing number of the virus, $R_0$.

 \section{Conclusion}
 \label{sec:ccl}

This paper proposes a set of novel indices adapted to target sensitivity analysis while being able to handle correlated inputs. The objective is to quantify the importance of inputs on the occurrence of a critical failure event of the system under study. The proposed indices are based on a cooperative Shapley procedure which aims at allocating the effects of the interaction and correlation equally between all the inputs in the same manner as the Shapley effects in global sensitivity analysis. Thus, a general class of distance-based indices is proposed, namely the $(\mathcal{D})$-target Shapley effects and some relevant properties are highlighted. Depending on the choice of the distance $\mathcal{D}$, well-known preexisting indices can be used as cost functions in the Shapley formulation. Therefore, these indices allow for the allocation, among the different inputs, of shares of several dispersion statistics (\eg mean absolute deviation for the $\ell^1$ case, variance for the $\ell^2$ case). \rev{These indices are} easily usable in practice, as they can be interpreted as percentages of the dispersion statistic, allocated to each input. This versatile procedure produces input importance measures according to a specific metric, driven by the choice of the distance.

In particular, the $(\ell^2)$-target Shapley effects (called target Shapley effects to simplify), which represents percentages of variance, have been studied more extensively and two dedicated estimation methods have been proposed. These particular indices have then been applied, analyzed and discussed through simple Gaussian toy-cases. Finally, two real-world use-cases have been studied: the modeling of a river flood and the ICU bed shortage during the COVID-19 pandemic. These indices are revealed to be able to detect influential inputs in the context of correlated inputs. For target sensitivity analysis, such a tool is valuable and can be used as a complement of more standard procedures. The clear advantage \rev{of this method} is that only one set of indices is \rev{needed} in order to produce easily interpretable and meaningful insights \rev{regarding the studied phenomenon}. Moreover, the proposed indices can be estimated in \rev{a} given-data \rev{context} which can be adapted to applications for which no computer model is available. However, the major limitations of the approach are primarily related to the target aspect of the analysis. Indeed, as soon as the event becomes increasingly rare, all the inputs tend to be influential and making a clear distinction between interactions and correlation effects becomes difficult.

To overcome these limitations, a first approach could be to improve the estimation strategies. The sampling-based method could benefit from a better sampling scheme, such as importance sampling, as described by \cite{RubinsteinKroese1981}, which could reduce the estimator's variance. Recent results from \cite{Sarazin_Derennes_Morio_AMM_2020} using copulas are also promising in the extent to which they show efficient estimations of the Shapley effects. Moreover, adapting recent results from \cite{spagnol_phd_2020}, with a link between the target Sobol' indices and the Squared Mutual Information \cite{sugiyama_machine_2012}, should allow for other possibilities of given-data estimation methods. Another method based on a random forest given-data procedure, explored by \cite{ElieDitCosaque_PHD_2020} in the context of quantile-oriented importance measure estimation, could also yield promising results if transposed to a reliability-oriented setting.

\rev{Even} if the Shapley attribution system is a solution when dealing with input statistical dependencies, it lacks a finer decomposition allowing to quantify the origin of each effect (\eg statistical dependence and interaction). \rev{Future} work \rev{could} use the recent developments in \cite{rabitti_shapley--owen_2019} in order to quantify interaction effects, \rev{by} transposition \rev{of these results} to the target sensitivity analysis setting.

Finally, it has been shown in \cite{soofi_framework_2000} that the Shapley attribution system is equivalent to a maximum entropy distribution (\eg uniform) over all possible orderings of inputs (the Shapley weights). Developments towards other forms of data-driven allocation systems could also open a path for further improvements.
\section*{Acknowledgments}

We are grateful to the three anonymous referees as well as Jérôme Morio, François Bachoc and Julien Demange-Chryst for their helpful remarks. We also thank S\'ebastien Da Veiga, Cl\'ementine Prieur and Fabrice Gamboa for interesting discussions and for having provided the dataset on the COVID-19 model. Finally, we would like to thank Victoria Stanford for her help in  proofreading this work.

\bibliographystyle{elsarticle-num-names}
\bibliography{bib/references, bib/biblio_VincentChabridon} 

\begin{thebibliography}{71}
\expandafter\ifx\csname natexlab\endcsname\relax\def\natexlab#1{#1}\fi
\providecommand{\url}[1]{\texttt{#1}}
\providecommand{\href}[2]{#2}
\providecommand{\path}[1]{#1}
\providecommand{\DOIprefix}{doi:}
\providecommand{\ArXivprefix}{arXiv:}
\providecommand{\URLprefix}{URL: }
\providecommand{\Pubmedprefix}{pmid:}
\providecommand{\doi}[1]{\href{http://dx.doi.org/#1}{\path{#1}}}
\providecommand{\Pubmed}[1]{\href{pmid:#1}{\path{#1}}}
\providecommand{\bibinfo}[2]{#2}
\ifx\xfnm\relax \def\xfnm[#1]{\unskip,\space#1}\fi
\bibitem[{De~Rocquigny et~al.(2008)De~Rocquigny, Devictor, and
  Tarantola}]{DeRocquigny_Devictor_book_2008}
\bibinfo{author}{E.~De~Rocquigny}, \bibinfo{author}{N.~Devictor},
  \bibinfo{author}{S.~Tarantola}, \bibinfo{title}{{Uncertainty in industrial
  practice: a guide to quantitative uncertainty management}},
  \bibinfo{publisher}{Wiley}, \bibinfo{year}{2008}.
\bibitem[{Beven(2008)}]{bev08}
\bibinfo{author}{K.~Beven}, \bibinfo{title}{Environmental Modelling: An
  Uncertain Future?}, \bibinfo{publisher}{CRC Press}, \bibinfo{year}{2008}.
\bibitem[{Pianosi et~al.(2016)Pianosi, Beven, Freer, Hall, Rougier, Stephenson,
  and Wagener}]{piabev16}
\bibinfo{author}{F.~Pianosi}, \bibinfo{author}{K.~Beven},
  \bibinfo{author}{J.~Freer}, \bibinfo{author}{J.~Hall},
  \bibinfo{author}{J.~Rougier}, \bibinfo{author}{D.~Stephenson},
  \bibinfo{author}{T.~Wagener},
\newblock \bibinfo{title}{Sensitivity analysis of environmental models: A
  systematic review with practical workflow},
\newblock \bibinfo{journal}{Environmental Modelling \& Software}
  \bibinfo{volume}{79} (\bibinfo{year}{2016}) \bibinfo{pages}{214--232}.
\bibitem[{Razavi et~al.(2021)Razavi, Jakeman, Saltelli, Prieur, Iooss,
  Borgonovo, Plischke, Lo~Piano, Iwanaga, Becker, Tarantola, Guillaume,
  Jakeman, Gupta, Melillo, Rabiti, Chabridon, Duan, Sun, Smith, Sheikholeslami,
  Hosseini, Asadzadeh, Puy, Kucherenko, and Maier}]{razjak20}
\bibinfo{author}{S.~Razavi}, \bibinfo{author}{A.~Jakeman},
  \bibinfo{author}{A.~Saltelli}, \bibinfo{author}{C.~Prieur},
  \bibinfo{author}{B.~Iooss}, \bibinfo{author}{E.~Borgonovo},
  \bibinfo{author}{E.~Plischke}, \bibinfo{author}{S.~Lo~Piano},
  \bibinfo{author}{T.~Iwanaga}, \bibinfo{author}{W.~Becker},
  \bibinfo{author}{S.~Tarantola}, \bibinfo{author}{J.~Guillaume},
  \bibinfo{author}{J.~Jakeman}, \bibinfo{author}{H.~Gupta},
  \bibinfo{author}{N.~Melillo}, \bibinfo{author}{G.~Rabiti},
  \bibinfo{author}{V.~Chabridon}, \bibinfo{author}{Q.~Duan},
  \bibinfo{author}{X.~Sun}, \bibinfo{author}{S.~Smith},
  \bibinfo{author}{R.~Sheikholeslami}, \bibinfo{author}{N.~Hosseini},
  \bibinfo{author}{M.~Asadzadeh}, \bibinfo{author}{A.~Puy},
  \bibinfo{author}{S.~Kucherenko}, \bibinfo{author}{H.~Maier},
\newblock \bibinfo{title}{The future of sensitivity analysis: An essential
  discipline for systems modelling and policy making},
\newblock \bibinfo{journal}{Environmental Modelling and Software}
  \bibinfo{volume}{137} (\bibinfo{year}{2021}) \bibinfo{pages}{104954}.
\bibitem[{Saltelli et~al.(2008)Saltelli, Ratto, Andres, Campolongo, Cariboni,
  Gatelli, Saisana, and Tarantola}]{Saltelli_THE_PRIMER_book}
\bibinfo{author}{A.~Saltelli}, \bibinfo{author}{M.~Ratto},
  \bibinfo{author}{T.~Andres}, \bibinfo{author}{F.~Campolongo},
  \bibinfo{author}{J.~Cariboni}, \bibinfo{author}{D.~Gatelli},
  \bibinfo{author}{M.~Saisana}, \bibinfo{author}{S.~Tarantola},
  \bibinfo{title}{{Global Sensitivity Analysis. The Primer}},
  \bibinfo{publisher}{Wiley}, \bibinfo{year}{2008}.
\bibitem[{Iooss and Lema\^{\i}tre(2015)}]{ioolem15}
\bibinfo{author}{B.~Iooss}, \bibinfo{author}{P.~Lema\^{\i}tre},
\newblock \bibinfo{title}{A review on global sensitivity analysis methods},
\newblock in: \bibinfo{editor}{C.~Meloni}, \bibinfo{editor}{G.~Dellino} (Eds.),
  \bibinfo{booktitle}{Uncertainty management in Simulation-Optimization of
  Complex Systems: Algorithms and Applications}, \bibinfo{publisher}{Springer},
  \bibinfo{year}{2015}, pp. \bibinfo{pages}{101--122}.
\bibitem[{Lemaire et~al.(2009)Lemaire, Chateauneuf, and Mitteau}]{Lemaire2009}
\bibinfo{author}{M.~Lemaire}, \bibinfo{author}{A.~Chateauneuf},
  \bibinfo{author}{J.-C. Mitteau}, \bibinfo{title}{{S}tructural {R}eliability},
  \bibinfo{publisher}{ISTE Ltd \& John Wiley \& Sons, Inc.},
  \bibinfo{year}{2009}.
\bibitem[{Richet and Bacchi(2019)}]{ricvit19}
\bibinfo{author}{Y.~Richet}, \bibinfo{author}{V.~Bacchi},
\newblock \bibinfo{title}{Inversion algorithm for civil flood defense
  optimization: Application to two-dimensional numerical model of the garonne
  river in france},
\newblock \bibinfo{journal}{Frontiers in Environmental Science}
  \bibinfo{volume}{7} (\bibinfo{year}{2019}) \bibinfo{pages}{160}.
\bibitem[{Rockafellar and Royset(2015)}]{Rockafellar_Royset_JRUES_2015}
\bibinfo{author}{R.~T. Rockafellar}, \bibinfo{author}{J.~O. Royset},
\newblock \bibinfo{title}{{Engineering Decisions under Risk Averseness}},
\newblock \bibinfo{journal}{ASCE-ASME Journal of Risk and Uncertainty in
  Engineering Systems, Part A: Civil Engineering} \bibinfo{volume}{1}
  (\bibinfo{year}{2015}) \bibinfo{pages}{1--12}.
\bibitem[{Morio and Balesdent(2015)}]{MorioBalesdent2015}
\bibinfo{author}{J.~Morio}, \bibinfo{author}{M.~Balesdent},
  \bibinfo{title}{{Estimation of Rare Event Probabilities in Complex Aerospace
  and Other Systems: A Practical Approach}}, \bibinfo{publisher}{Woodhead
  Publishing, Elsevier}, \bibinfo{year}{2015}.
\bibitem[{Wu(1994)}]{Wu1994}
\bibinfo{author}{Y.-T. Wu},
\newblock \bibinfo{title}{{Computational Methods for Efficient Structural
  Reliability and Reliability Sensitivity Analysis}},
\newblock \bibinfo{journal}{AIAA Journal} \bibinfo{volume}{32}
  (\bibinfo{year}{1994}) \bibinfo{pages}{1717--1723}.
\bibitem[{Song et~al.(2009)Song, Lu, and Qiao}]{Song2009}
\bibinfo{author}{S.~Song}, \bibinfo{author}{Z.~Lu}, \bibinfo{author}{H.~Qiao},
\newblock \bibinfo{title}{{Subset simulation for structural reliability
  sensitivity analysis}},
\newblock \bibinfo{journal}{Reliability Engineering and System Safety}
  \bibinfo{volume}{94} (\bibinfo{year}{2009}) \bibinfo{pages}{658--665}.
\bibitem[{Wei et~al.(2012)Wei, Lu, Hao, Feng, and Wang}]{Wei_Lu_CPC_2012}
\bibinfo{author}{P.~Wei}, \bibinfo{author}{Z.~Lu}, \bibinfo{author}{W.~Hao},
  \bibinfo{author}{J.~Feng}, \bibinfo{author}{B.~Wang},
\newblock \bibinfo{title}{{Efficient sampling methods for global reliability
  sensitivity analysis}},
\newblock \bibinfo{journal}{Computer Physics Communications}
  \bibinfo{volume}{183} (\bibinfo{year}{2012}) \bibinfo{pages}{1728--1743}.
\bibitem[{Chabridon(2018)}]{chabridon_phd_2018}
\bibinfo{author}{V.~Chabridon}, \bibinfo{title}{{Reliability-oriented
  sensitivity analysis under probabilistic model uncertainty -- Application to
  aerospace systems}}, Ph.D. thesis, Universit\'e Clermont Auvergne,
  \bibinfo{year}{2018}.
\bibitem[{Perrin and Defaux(2019)}]{Perrin_Defaux_JSC_2019}
\bibinfo{author}{G.~Perrin}, \bibinfo{author}{G.~Defaux},
\newblock \bibinfo{title}{{Efficient Evaluation of Reliability-Oriented
  Sensitivity Indices}},
\newblock \bibinfo{journal}{Journal of Scientific Computing}
  \bibinfo{volume}{79} (\bibinfo{year}{2019}) \bibinfo{pages}{1433--1455}.
\bibitem[{Sobol(1993)}]{Sobol_MMCE_1993}
\bibinfo{author}{I.~M. Sobol},
\newblock \bibinfo{title}{{Sensitivity estimates for nonlinear mathematical
  models}},
\newblock \bibinfo{journal}{Mathematical Modelling and Computational
  Experiments} \bibinfo{volume}{1} (\bibinfo{year}{1993})
  \bibinfo{pages}{407--414}.
\bibitem[{Borgonovo(2007)}]{borgonovo_new_2007}
\bibinfo{author}{E.~Borgonovo},
\newblock \bibinfo{title}{A new uncertainty importance measure},
\newblock \bibinfo{journal}{Reliability Engineering \& System Safety}
  \bibinfo{volume}{92} (\bibinfo{year}{2007}) \bibinfo{pages}{771--784}.
\bibitem[{Raguet and Marrel(2018)}]{raguet_target_2018}
\bibinfo{author}{H.~Raguet}, \bibinfo{author}{A.~Marrel},
\newblock \bibinfo{title}{Target and conditional sensitivity analysis with
  emphasis on dependence measures},
\newblock \bibinfo{journal}{Working paper}  (\bibinfo{year}{2018}). \URLprefix
  \url{https://arxiv.org/abs/1801.10047}.
\bibitem[{Li et~al.(2012)Li, Lu, Jun, and Bintuan}]{Li_Lu_SS_2012}
\bibinfo{author}{L.~Li}, \bibinfo{author}{Z.~Lu}, \bibinfo{author}{F.~Jun},
  \bibinfo{author}{W.~Bintuan},
\newblock \bibinfo{title}{{Moment-independent importance measure of basic
  variable and its state dependent parameter solution}},
\newblock \bibinfo{journal}{Structural Safety} \bibinfo{volume}{38}
  (\bibinfo{year}{2012}) \bibinfo{pages}{40--47}.
\bibitem[{Marrel and Chabridon(2021)}]{marrel_statistical_2020}
\bibinfo{author}{A.~Marrel}, \bibinfo{author}{V.~Chabridon},
\newblock \bibinfo{title}{{Statistical developments for target and conditional
  sensitivity analysis: application on safety studies for nuclear reactor}},
\newblock \bibinfo{journal}{Reliability Engineering and System Safety, in
  press}  (\bibinfo{year}{2021}). \URLprefix
  \url{https://hal.archives-ouvertes.fr/hal-02541142v2}.
\bibitem[{Hoeffding(1948)}]{Hoeffding_1948}
\bibinfo{author}{W.~Hoeffding},
\newblock \bibinfo{title}{{A class of statistics with asymptotically normal
  distribution}},
\newblock \bibinfo{journal}{The Annals of Mathematical Statistics}
  \bibinfo{volume}{19} (\bibinfo{year}{1948}) \bibinfo{pages}{293--325}.
\bibitem[{Jacques et~al.(2006)Jacques, Lavergne, and
  Devictor}]{jacques_sensitivity_2006}
\bibinfo{author}{J.~Jacques}, \bibinfo{author}{C.~Lavergne},
  \bibinfo{author}{N.~Devictor},
\newblock \bibinfo{title}{Sensitivity analysis in presence of model uncertainty
  and correlated inputs},
\newblock \bibinfo{journal}{Reliability Engineering \& System Safety}
  \bibinfo{volume}{91} (\bibinfo{year}{2006}) \bibinfo{pages}{1126--1134}.
\bibitem[{Li et~al.(2010)Li, Rabitz, Yelvington, Oluwole, Bacon, Kolb, and
  Schoendorf}]{lirab10}
\bibinfo{author}{G.~Li}, \bibinfo{author}{H.~Rabitz},
  \bibinfo{author}{P.~Yelvington}, \bibinfo{author}{O.~Oluwole},
  \bibinfo{author}{F.~Bacon}, \bibinfo{author}{C.~Kolb},
  \bibinfo{author}{J.~Schoendorf},
\newblock \bibinfo{title}{Global sensitivity analysis for systems with
  independent and/or correlated inputs},
\newblock \bibinfo{journal}{Journal of Physical Chemistry}
  \bibinfo{volume}{114} (\bibinfo{year}{2010}) \bibinfo{pages}{6022--6032}.
\bibitem[{Chastaing et~al.(2012)Chastaing, Gamboa, and
  Prieur}]{chastaing_generalized_2012}
\bibinfo{author}{G.~Chastaing}, \bibinfo{author}{F.~Gamboa},
  \bibinfo{author}{C.~Prieur},
\newblock \bibinfo{title}{Generalized {H}oeffding-{S}obol decomposition for
  dependent variables - {A}pplication to sensitivity analysis},
\newblock \bibinfo{journal}{Electronic Journal of Statistics}
  \bibinfo{volume}{6} (\bibinfo{year}{2012}) \bibinfo{pages}{2420--2448}.
\bibitem[{Xu and Gertner(2008)}]{xu_uncertainty_2008}
\bibinfo{author}{C.~Xu}, \bibinfo{author}{G.~Z. Gertner},
\newblock \bibinfo{title}{Uncertainty and sensitivity analysis for models with
  correlated parameters},
\newblock \bibinfo{journal}{Reliability Engineering \& System Safety}
  \bibinfo{volume}{93} (\bibinfo{year}{2008}) \bibinfo{pages}{1563--1573}.
\bibitem[{Mara and Tarantola(2012)}]{martar12}
\bibinfo{author}{T.~Mara}, \bibinfo{author}{S.~Tarantola},
\newblock \bibinfo{title}{Variance-based sensitivity indices for models with
  dependent inputs},
\newblock \bibinfo{journal}{Reliability Engineering \& System Safety}
  \bibinfo{volume}{107} (\bibinfo{year}{2012}) \bibinfo{pages}{115--121}.
\bibitem[{Mara et~al.(2015)Mara, Tarantola, and Annoni}]{martar15}
\bibinfo{author}{T.~Mara}, \bibinfo{author}{S.~Tarantola},
  \bibinfo{author}{P.~Annoni},
\newblock \bibinfo{title}{Non-parametric methods for global sensitivity
  analysis of model output with dependent inputs},
\newblock \bibinfo{journal}{Environmental Modeling \& Software}
  \bibinfo{volume}{72} (\bibinfo{year}{2015}) \bibinfo{pages}{173--183}.
\bibitem[{Benoumechiara and Elie-Dit-Cosaque(2019)}]{beneli19}
\bibinfo{author}{N.~Benoumechiara}, \bibinfo{author}{K.~Elie-Dit-Cosaque},
\newblock \bibinfo{title}{Shapley effects for sensitivity analysis with
  dependent inputs: bootstrap and kriging-based algorithms},
\newblock \bibinfo{journal}{ESAIM: Proceedings and Surveys}
  \bibinfo{volume}{65} (\bibinfo{year}{2019}) \bibinfo{pages}{266--293}.
\bibitem[{Do and Razavi(2020)}]{doraz20}
\bibinfo{author}{N.~Do}, \bibinfo{author}{S.~Razavi},
\newblock \bibinfo{title}{Correlation effects? {A} major but often neglected
  component in sensitivity and uncertainty analysis},
\newblock \bibinfo{journal}{Water Resources Research} \bibinfo{volume}{56}
  (\bibinfo{year}{2020}) \bibinfo{pages}{e2019WR025436}.
\bibitem[{Shapley(1953)}]{Shapley_1953}
\bibinfo{author}{L.~S. Shapley},
\newblock \bibinfo{title}{{A value for n-person games}},
\newblock in: \bibinfo{editor}{H.~Kuhn}, \bibinfo{editor}{A.~W. Tucker} (Eds.),
  \bibinfo{booktitle}{{Contributions to the Theory of Games, Volume II}},
  Annals of Mathematics Studies, \bibinfo{publisher}{Princeton University
  Press}, \bibinfo{address}{Princeton, NJ}, \bibinfo{year}{1953}, pp.
  \bibinfo{pages}{307--317}.
\bibitem[{Osborne and Rubinstein(1994)}]{osborne_course_1994}
\bibinfo{author}{M.~Osborne}, \bibinfo{author}{A.~Rubinstein},
  \bibinfo{title}{{A Course in Game Theory}}, \bibinfo{publisher}{{MIT} Press},
  \bibinfo{year}{1994}.
\bibitem[{Owen(2014)}]{Owen_ASAJUQ_2014}
\bibinfo{author}{A.~B. Owen},
\newblock \bibinfo{title}{{Sobol' indices and Shapley value}},
\newblock \bibinfo{journal}{SIAM/ASA Journal of Uncertainty Quantification}
  \bibinfo{volume}{2} (\bibinfo{year}{2014}) \bibinfo{pages}{245--251}.
\bibitem[{Owen and Prieur(2017)}]{Owen_Prieur_ASAJUQ_2017}
\bibinfo{author}{A.~B. Owen}, \bibinfo{author}{C.~Prieur},
\newblock \bibinfo{title}{{On Shapley value for measuring importance of
  dependent inputs}},
\newblock \bibinfo{journal}{SIAM/ASA Journal of Uncertainty Quantification}
  \bibinfo{volume}{5} (\bibinfo{year}{2017}) \bibinfo{pages}{986--1002}.
\bibitem[{Iooss and Prieur(2019)}]{iooss_shapley_2019}
\bibinfo{author}{B.~Iooss}, \bibinfo{author}{C.~Prieur},
\newblock \bibinfo{title}{{Shapley effects for Sensitivity Analysis with
  correlated inputs : Comparisons with Sobol' Indices, Numerical Estimation and
  Applications}},
\newblock \bibinfo{journal}{International Journal for Uncertainty
  Quantification} \bibinfo{volume}{9} (\bibinfo{year}{2019})
  \bibinfo{pages}{493--514}.
\bibitem[{Spagnol(2020)}]{spagnol_phd_2020}
\bibinfo{author}{A.~Spagnol}, \bibinfo{title}{Kernel-based sensitivity indices
  for high-dimensional optimization problems}, Ph.D. thesis, Ecole des Mines de
  Saint-Etienne, \bibinfo{year}{2020}.
\bibitem[{Broto et~al.(2020)Broto, Bachoc, and Depecker}]{broto_variance_2020}
\bibinfo{author}{B.~Broto}, \bibinfo{author}{F.~Bachoc},
  \bibinfo{author}{M.~Depecker},
\newblock \bibinfo{title}{Variance reduction for estimation of shapley effects
  and adaptation to unknown input distribution},
\newblock \bibinfo{journal}{{SIAM}/{ASA} Journal on Uncertainty Quantification}
  \bibinfo{volume}{8} (\bibinfo{year}{2020}) \bibinfo{pages}{693--716}.
\bibitem[{Christensen(1990)}]{chr90}
\bibinfo{author}{R.~Christensen}, \bibinfo{title}{Linear models for
  multivariate, time series and spatial data},
  \bibinfo{publisher}{Springer-Verlag}, \bibinfo{year}{1990}.
\bibitem[{Hastie et~al.(2002)Hastie, Tibshirani, and Friedman}]{hastib02}
\bibinfo{author}{T.~Hastie}, \bibinfo{author}{R.~Tibshirani},
  \bibinfo{author}{J.~Friedman}, \bibinfo{title}{The elements of statistical
  learning}, \bibinfo{publisher}{Springer}, \bibinfo{year}{2002}.
\bibitem[{Helton et~al.(2006)Helton, Johnson, Salaberry, and
  Storlie}]{heljoh06}
\bibinfo{author}{J.~Helton}, \bibinfo{author}{J.~Johnson},
  \bibinfo{author}{C.~Salaberry}, \bibinfo{author}{C.~Storlie},
\newblock \bibinfo{title}{Survey of sampling-based methods for uncertainty and
  sensitivity analysis},
\newblock \bibinfo{journal}{Reliability Engineering \& System Safety}
  \bibinfo{volume}{91} (\bibinfo{year}{2006}) \bibinfo{pages}{1175--1209}.
\bibitem[{Johnson and LeBreton(2004)}]{johleb04}
\bibinfo{author}{J.~Johnson}, \bibinfo{author}{J.~LeBreton},
\newblock \bibinfo{title}{History and use of relative importance indices in
  organizational research},
\newblock \bibinfo{journal}{Organizational Research Methods}
  \bibinfo{volume}{7} (\bibinfo{year}{2004}) \bibinfo{pages}{238--257}.
\bibitem[{Clouvel(2019)}]{clo19}
\bibinfo{author}{L.~Clouvel}, \bibinfo{title}{Uncertainty quantification of the
  fast flux calculation for a {PWR} vessel}, Ph.D. thesis, Universit\'e
  Paris-Saclay, \bibinfo{year}{2019}.
\bibitem[{Lindeman et~al.(1980)Lindeman, Merenda, and Gold}]{linmer80}
\bibinfo{author}{R.~H. Lindeman}, \bibinfo{author}{P.~F. Merenda},
  \bibinfo{author}{R.~Z. Gold}, \bibinfo{title}{Introduction to bivariate and
  multivariate analysis}, \bibinfo{publisher}{Scott Foresman and Company},
  \bibinfo{address}{Glenview, IL}, \bibinfo{year}{1980}.
\bibitem[{Gr{\"o}mping(2006)}]{gro06}
\bibinfo{author}{U.~Gr{\"o}mping},
\newblock \bibinfo{title}{Relative importance for linear regression in
  \textsf{R}: the {P}ackage \textsf{relaimpo}},
\newblock \bibinfo{journal}{Journal of Statistical Software}
  \bibinfo{volume}{17} (\bibinfo{year}{2006}) \bibinfo{pages}{1--27}.
\bibitem[{Nossent et~al.(2011)Nossent, Elsen, and Bauwens}]{nosels11}
\bibinfo{author}{J.~Nossent}, \bibinfo{author}{P.~Elsen},
  \bibinfo{author}{W.~Bauwens},
\newblock \bibinfo{title}{Sobol’ sensitivity analysis of a complex
  environmental model},
\newblock \bibinfo{journal}{Environmental Modelling \& Software}
  \bibinfo{volume}{26} (\bibinfo{year}{2011}) \bibinfo{pages}{1515 -- 1525}.
\bibitem[{Song et~al.(2016)Song, Nelson, and Staum}]{song_shapley_2016}
\bibinfo{author}{E.~Song}, \bibinfo{author}{B.~Nelson},
  \bibinfo{author}{J.~Staum},
\newblock \bibinfo{title}{Shapley effects for global sensitivity analysis:
  Theory and computation},
\newblock \bibinfo{journal}{{SIAM}/{ASA} Journal on Uncertainty Quantification}
  \bibinfo{volume}{4} (\bibinfo{year}{2016}) \bibinfo{pages}{1060--1083}.
\bibitem[{Brandenburger(2007)}]{brandenburger_cooperative_2007}
\bibinfo{author}{A.~Brandenburger}, \bibinfo{title}{Cooperative {Game}
  {Theory}: {Characteristic} {Functions}, {Allocations}, {Marginal}
  {Contribution}}, \bibinfo{year}{2007}. \URLprefix
  \url{https://web.archive.org/web/20170829122711/http://www.uib.cat/depart/deeweb/pdi/hdeelbm0/arxius_decisions_and_games/cooperative_game_theory-brandenburger.pdf}.
\bibitem[{Derennes et~al.(2021)Derennes, Morio, and
  Simatos}]{Derennes_Morio_Simatos_MACS_2021}
\bibinfo{author}{P.~Derennes}, \bibinfo{author}{J.~Morio},
  \bibinfo{author}{F.~Simatos},
\newblock \bibinfo{title}{{Simultaneous estimation of complementary moment
  independent and reliability-oriented sensitivity measures}},
\newblock \bibinfo{journal}{Mathematics and Computers in Simulation}
  \bibinfo{volume}{182} (\bibinfo{year}{2021}) \bibinfo{pages}{721--737}.
\bibitem[{Morio(2012)}]{Morio_SMPT_2012}
\bibinfo{author}{J.~Morio},
\newblock \bibinfo{title}{{Extreme quantile estimation with nonparametric
  adaptive importance sampling}},
\newblock \bibinfo{journal}{Simulation Modelling Practice and Theory}
  \bibinfo{volume}{27} (\bibinfo{year}{2012}) \bibinfo{pages}{76--89}.
\bibitem[{Chabridon et~al.(2020)Chabridon, Balesdent, Perrin, Morio, Bourinet,
  and Gayton}]{Chabridon_Chapter_2020}
\bibinfo{author}{V.~Chabridon}, \bibinfo{author}{M.~Balesdent},
  \bibinfo{author}{G.~Perrin}, \bibinfo{author}{J.~Morio},
  \bibinfo{author}{J.-M. Bourinet}, \bibinfo{author}{N.~Gayton},
  \bibinfo{title}{{Mechanical Engineering Under Uncertainties}},
  \bibinfo{publisher}{Wiley - ISTE Ltd}, \bibinfo{year}{2020}, pp.
  \bibinfo{pages}{1--43}.
\bibitem[{Cui et~al.(2010)Cui, Lu, and Zhao}]{Cui_Lu_SCTS_2010}
\bibinfo{author}{L.~Cui}, \bibinfo{author}{Z.~Lu}, \bibinfo{author}{X.~Zhao},
\newblock \bibinfo{title}{{Moment-independent importance measure of basic
  random variable and its probability density evolution solution}},
\newblock \bibinfo{journal}{Science China Technical Sciences}
  \bibinfo{volume}{53} (\bibinfo{year}{2010}) \bibinfo{pages}{1138--1145}.
\bibitem[{Fort et~al.(2016)Fort, Klein, and Rachdi}]{FortKleinRachdi2016}
\bibinfo{author}{J.-C. Fort}, \bibinfo{author}{T.~Klein},
  \bibinfo{author}{N.~Rachdi},
\newblock \bibinfo{title}{{New sensitivity analysis subordinated to a
  contrast}},
\newblock \bibinfo{journal}{Communications in Statistics - Theory and Methods}
  \bibinfo{volume}{45} (\bibinfo{year}{2016}) \bibinfo{pages}{4349--4364}.
\bibitem[{Browne et~al.(2017)Browne, Fort, Iooss, and
  Le~Gratiet}]{BrowneFortIoossLeGratiet_sensi2017}
\bibinfo{author}{T.~Browne}, \bibinfo{author}{J.-C. Fort},
  \bibinfo{author}{B.~Iooss}, \bibinfo{author}{L.~Le~Gratiet},
\newblock \bibinfo{title}{{Estimate of quantile-oriented sensitivity indices}},
\newblock \bibinfo{journal}{HAL, hal-01450891, version 1}
  (\bibinfo{year}{2017}).
\bibitem[{Maume-Deschamps and Niang(2018)}]{MaumeDeschampsNiang_2018}
\bibinfo{author}{V.~Maume-Deschamps}, \bibinfo{author}{I.~Niang},
\newblock \bibinfo{title}{{Estimation of quantile oriented sensitivity
  indices}},
\newblock \bibinfo{journal}{Statistics and Probability Letters}
  \bibinfo{volume}{134} (\bibinfo{year}{2018}) \bibinfo{pages}{122--127}.
\bibitem[{Kucherenko et~al.(2019)Kucherenko, Song, and
  Wang}]{Kucherenko_Song_Wang_RESS_2019}
\bibinfo{author}{S.~Kucherenko}, \bibinfo{author}{S.~Song},
  \bibinfo{author}{L.~Wang},
\newblock \bibinfo{title}{{Quantile based global sensitivity measures}},
\newblock \bibinfo{journal}{Reliability Engineering and System Safety}
  \bibinfo{volume}{185} (\bibinfo{year}{2019}) \bibinfo{pages}{35--48}.
\bibitem[{Li et~al.(2016)Li, Lu, and Chen}]{Li_Lu_AST_2016}
\bibinfo{author}{L.~Li}, \bibinfo{author}{Z.~Lu}, \bibinfo{author}{C.~Chen},
\newblock \bibinfo{title}{{Moment-independent importance measure of correlated
  input variable and its state dependent parameter solution}},
\newblock \bibinfo{journal}{Aerospace Science and Technology}
  \bibinfo{volume}{48} (\bibinfo{year}{2016}) \bibinfo{pages}{281--290}.
\bibitem[{Iooss et~al.(2021)Iooss, {Da~Veiga}, Janon, and Pujol}]{ioodav20}
\bibinfo{author}{B.~Iooss}, \bibinfo{author}{S.~{Da~Veiga}},
  \bibinfo{author}{A.~Janon}, \bibinfo{author}{G.~Pujol},
  \bibinfo{title}{sensitivity: Global Sensitivity Analysis of Model Outputs},
  \bibinfo{year}{2021}. \URLprefix
  \url{https://CRAN.R-project.org/package=sensitivity}, \bibinfo{note}{{R}
  package version 1.25.0}.
\bibitem[{Lemaitre(2014)}]{lemaitre_analyse_2014}
\bibinfo{author}{P.~Lemaitre}, \bibinfo{title}{Analyse de sensibilité en
  fiabilité des structures}, Ph.D. thesis, Université de Bordeaux,
  \bibinfo{year}{2014}.
\bibitem[{Fox and Monette(1992)}]{foxmon92}
\bibinfo{author}{J.~Fox}, \bibinfo{author}{G.~Monette},
\newblock \bibinfo{title}{Generalized collinearity diagnostics},
\newblock \bibinfo{journal}{Journal of the American Statistical Association}
  \bibinfo{volume}{87} (\bibinfo{year}{1992}) \bibinfo{pages}{178--183}.
\bibitem[{Ortmann(2000)}]{ortmann_proportional_2000}
\bibinfo{author}{K.~M. Ortmann},
\newblock \bibinfo{title}{The proportional value for positive cooperative
  games},
\newblock \bibinfo{journal}{Mathematical Methods of Operations Research (ZOR)}
  \bibinfo{volume}{51} (\bibinfo{year}{2000}) \bibinfo{pages}{235--248}.
  \URLprefix \url{http://link.springer.com/10.1007/s001860050086}.
  \DOIprefix\doi{10.1007/s001860050086}.
\bibitem[{Feldman(2005)}]{feldman_relative_2005}
\bibinfo{author}{B.~E. Feldman},
\newblock \bibinfo{title}{Relative {Importance} and {Value}},
\newblock \bibinfo{journal}{SSRN Electronic Journal}  (\bibinfo{year}{2005}).
  \URLprefix \url{http://www.ssrn.com/abstract=2255827}.
  \DOIprefix\doi{10.2139/ssrn.2255827}.
\bibitem[{Schumann(2009)}]{schumann_algo_2009}
\bibinfo{author}{E.~Schumann},
\newblock \bibinfo{title}{Generating correlated uniform variates.}
  (\bibinfo{year}{2009}). \URLprefix
  \url{http://comisef.wikidot.com/tutorial:correlateduniformvariates}.
\bibitem[{Saltelli et~al.(2020)Saltelli, Bammer, Bruno, Charters, Fiore
  et~al.}]{salbam20}
\bibinfo{author}{A.~Saltelli}, \bibinfo{author}{G.~Bammer},
  \bibinfo{author}{I.~Bruno}, \bibinfo{author}{E.~Charters},
  \bibinfo{author}{M.~D. Fiore}, et~al.,
\newblock \bibinfo{title}{Five ways to ensure that models serve society: a
  manifesto (short comments)},
\newblock \bibinfo{journal}{Nature} \bibinfo{volume}{582}
  (\bibinfo{year}{2020}) \bibinfo{pages}{482--484}.
\bibitem[{Lu and Borgonovo(2020)}]{lubor20}
\bibinfo{author}{X.~Lu}, \bibinfo{author}{E.~Borgonovo},
\newblock \bibinfo{title}{Is time to intervention in the {COVID}-19 outbreak
  really important? {A} global sensitivity analysis approach},
\newblock \bibinfo{journal}{Preprint}  (\bibinfo{year}{2020}).
  \bibinfo{note}{ArXiv:2005.01833}.
\bibitem[{Da~Veiga et~al.(2021)Da~Veiga, Gamboa, Iooss, and
  Prieur}]{veiga_basics_2020}
\bibinfo{author}{S.~Da~Veiga}, \bibinfo{author}{F.~Gamboa},
  \bibinfo{author}{B.~Iooss}, \bibinfo{author}{C.~Prieur},
  \bibinfo{title}{Basics and trends in sensitivity analysis},
  \bibinfo{publisher}{SIAM, In press}, \bibinfo{year}{2021}.
\bibitem[{{Da Veiga}(2020)}]{dav20}
\bibinfo{author}{S.~{Da Veiga}},
\newblock \bibinfo{title}{Calibration and sensitivity analysis of a {COVID}-19
  epidemics model},
\newblock \bibinfo{organization}{Meeting AppliBUGS (Applications du Bayesian
  Unified Group of Statisticians)}, \bibinfo{year}{December 2020}. \URLprefix
  \url{genome.jouy.inra.fr/applibugs/Daveiga_AppliBUGSDec2020.pdf}.
\bibitem[{Rubinstein and Kroese(2008)}]{RubinsteinKroese1981}
\bibinfo{author}{R.~Y. Rubinstein}, \bibinfo{author}{D.~P. Kroese},
  \bibinfo{title}{{S}imulation and the {M}onte {C}arlo method},
  \bibinfo{edition}{{S}econd} ed., \bibinfo{publisher}{Wiley},
  \bibinfo{year}{2008}.
\bibitem[{Sarazin et~al.(2020)Sarazin, Derennes, and
  Morio}]{Sarazin_Derennes_Morio_AMM_2020}
\bibinfo{author}{G.~Sarazin}, \bibinfo{author}{P.~Derennes},
  \bibinfo{author}{J.~Morio},
\newblock \bibinfo{title}{{Estimation of high-order moment-independent
  importance measures for Shapley value analysis}},
\newblock \bibinfo{journal}{Applied Mathematical Modelling}
  \bibinfo{volume}{88} (\bibinfo{year}{2020}) \bibinfo{pages}{396--417}.
\bibitem[{Sugiyama(2012)}]{sugiyama_machine_2012}
\bibinfo{author}{M.~Sugiyama},
\newblock \bibinfo{title}{Machine learning with squared-loss mutual
  information},
\newblock \bibinfo{journal}{Entropy} \bibinfo{volume}{15}
  (\bibinfo{year}{2012}) \bibinfo{pages}{80--112}.
\bibitem[{Elie-Dit-Cosaque(2020)}]{ElieDitCosaque_PHD_2020}
\bibinfo{author}{K.~Elie-Dit-Cosaque}, \bibinfo{title}{{Développement de
  mesures d’incertitudes pour le risque de modèle dans des contextes
  incluant de la dépendance stochastique}}, Ph.D. thesis, Universit\'e Claude
  Bernard - Lyon 1, \bibinfo{year}{2020}.
\bibitem[{Rabitti and Borgonovo(2019)}]{rabitti_shapley--owen_2019}
\bibinfo{author}{G.~Rabitti}, \bibinfo{author}{E.~Borgonovo},
\newblock \bibinfo{title}{A {S}hapley--owen index for interaction
  quantification},
\newblock \bibinfo{journal}{{SIAM}/{ASA} Journal on Uncertainty Quantification}
  \bibinfo{volume}{7} (\bibinfo{year}{2019}) \bibinfo{pages}{1060--1075}.
\bibitem[{Soofi et~al.(2000)Soofi, Retzer, and
  Yasai-Ardekani}]{soofi_framework_2000}
\bibinfo{author}{E.~S. Soofi}, \bibinfo{author}{J.~J. Retzer},
  \bibinfo{author}{M.~Yasai-Ardekani},
\newblock \bibinfo{title}{A {Framework} for {Measuring} the {Importance} of
  {Variables} with {Applications} to {Management} {Research} and {Decision}
  {Models}*},
\newblock \bibinfo{journal}{Decision Sciences} \bibinfo{volume}{31}
  (\bibinfo{year}{2000}) \bibinfo{pages}{595 -- 625}.

\end{thebibliography}
\appendix
\renewcommand\thefigure{\thesection.\arabic{figure}} 
\section{ANOVA and Sobol' indices}
\setcounter{figure}{0}    
\label{app:sobol}



In the general non-linear case, as for the ANOVA of the linear model case (see Subsection \ref{sec:linear}), the idea is to find a general decomposition of the output variance. 
This can be done through the decomposition of a function with finite variance ($L^2$ mathematical property), called the \textit{Hoeffding decomposition} \cite{Hoeffding_1948}, which allows to rewrite $G(X)$ as a sum of centered components related to each possible subset of inputs. 
For example, in the case of a model with three inputs $X=(X_1, X_2, X_3)$, $G(X)$ can be decomposed into four components:
\begin{align}
\label{eq:hoef_dec}
G(X) &= G_{\emptyset} \tag*{(Mean behavior)}\\
& + G_{1}(X_1) + G_{2}(X_2) +G_{3}(X_3) \tag*{(First-order)}\\
&+  G_{\{1,2\}}(X_1, X_2) +  G_{\{1,3\}}(X_1, X_3) +  G_{\{2,3\}}(X_2, X_3) \tag*{(Second-order)}\\
&+G_{\{1,2,3\}}(X)\;. \tag*{(Third-order)}
\end{align}

Moreover, if the inputs are assumed to be independent, each term is orthogonal to one another and writes
\begin{equation}
\label{eq:hoef_sol}
G_{A}(x_A) = \sum_{B \subset A} (-1)^{|A| - |B|} \mathbb{E}\left[ G(X) | X_B = x_B\right]
\end{equation}
where $A \in \mathcal{P}_d$ is a subset of indices and $\mathcal{P}_d$ the set of all possible subsets of $\{1, \dots, d\}$, $|A|$ is the cardinal of $A$ and $X_A$ denotes the subset of inputs, selected by the indices in $A$ ($X_A = (X_i)_{i \in A}$).
Then, the Hoeffding decomposition is unique and leads to a variance decomposition called ``functional ANOVA'':
\begin{equation}\label{eq:FANOVA}
\mathbb{V}[G(X)] = \sum_{A \in \mathcal{P}_d,A\neq 0} \mathbb{V}[G_{A}(x_A)]\;.
\end{equation}
\noindent
This leads to the definition of the Sobol' indices \cite{Sobol_MMCE_1993}:
\begin{equation}
\label{eq:sobol}
S_A = \frac{\mathbb{V}[G_A(X_A)]}{\mathbb{V}[G(X)]} = \frac{\sum_{B \subset A} (-1)^{|A| - |B|} \mathbb{V}\left(\mathbb{E}\bigl[ G(X) \bigm| X_B\bigr]\right)}{\mathbb{V}[G(X)]}\;.
\end{equation}
The sum of the Sobol' indices over all subset on inputs $A \in \mathcal{P}_d$ being equal to one, they can be directly interpreted as the percentage of the output variance due to each subset of input \cite{Sobol_MMCE_1993,Saltelli_THE_PRIMER_book}.
The Sobol' indices of higher orders than one can be interpreted as a means of quantifying the share of variance due to the interaction effects induced by the structure of the model $G(\cdot)$ between the selected subset of inputs.

\begin{figure}[!ht]
\centering
\includegraphics[width=0.4\textwidth]{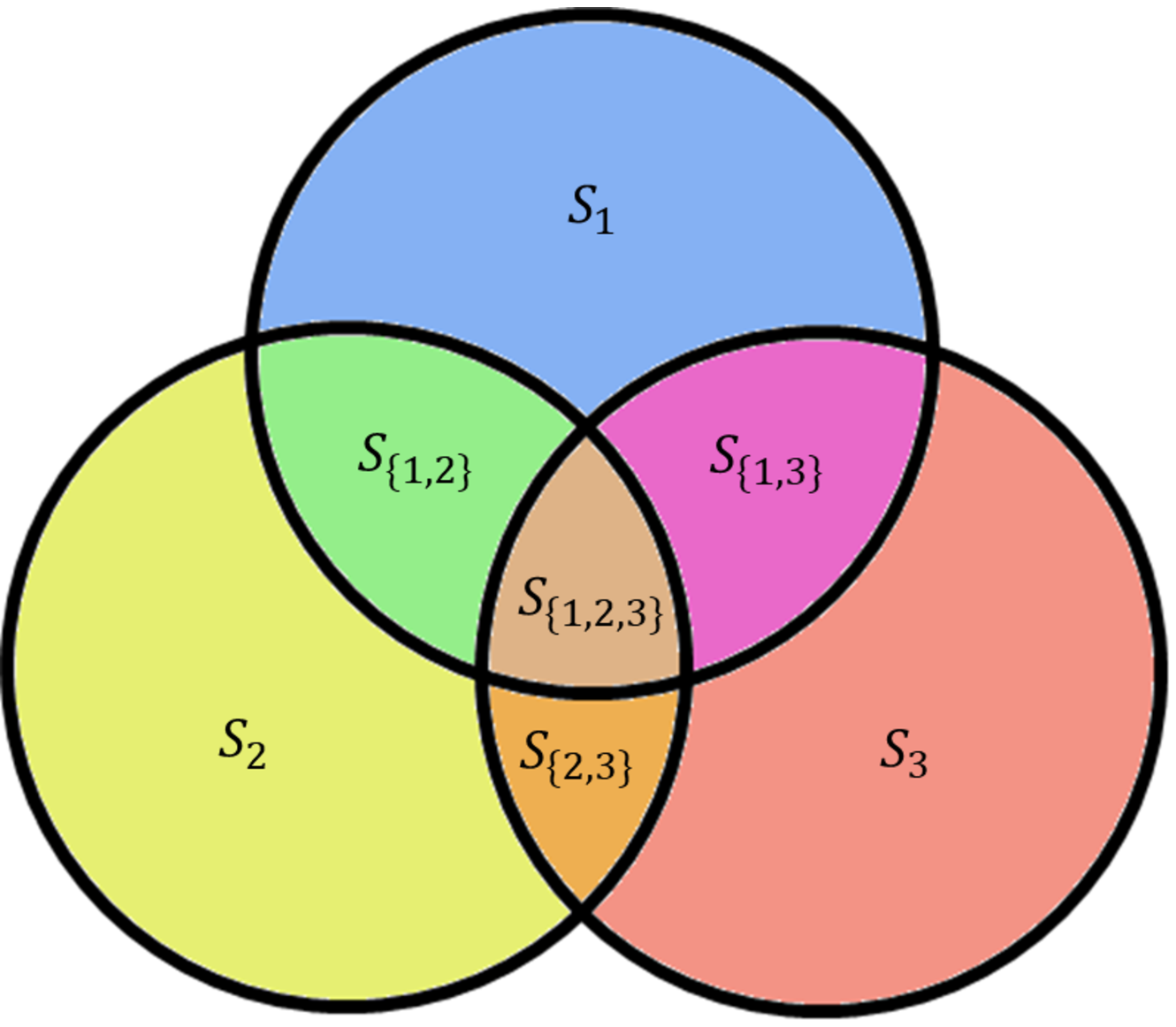}
\includegraphics[width=0.49\textwidth]{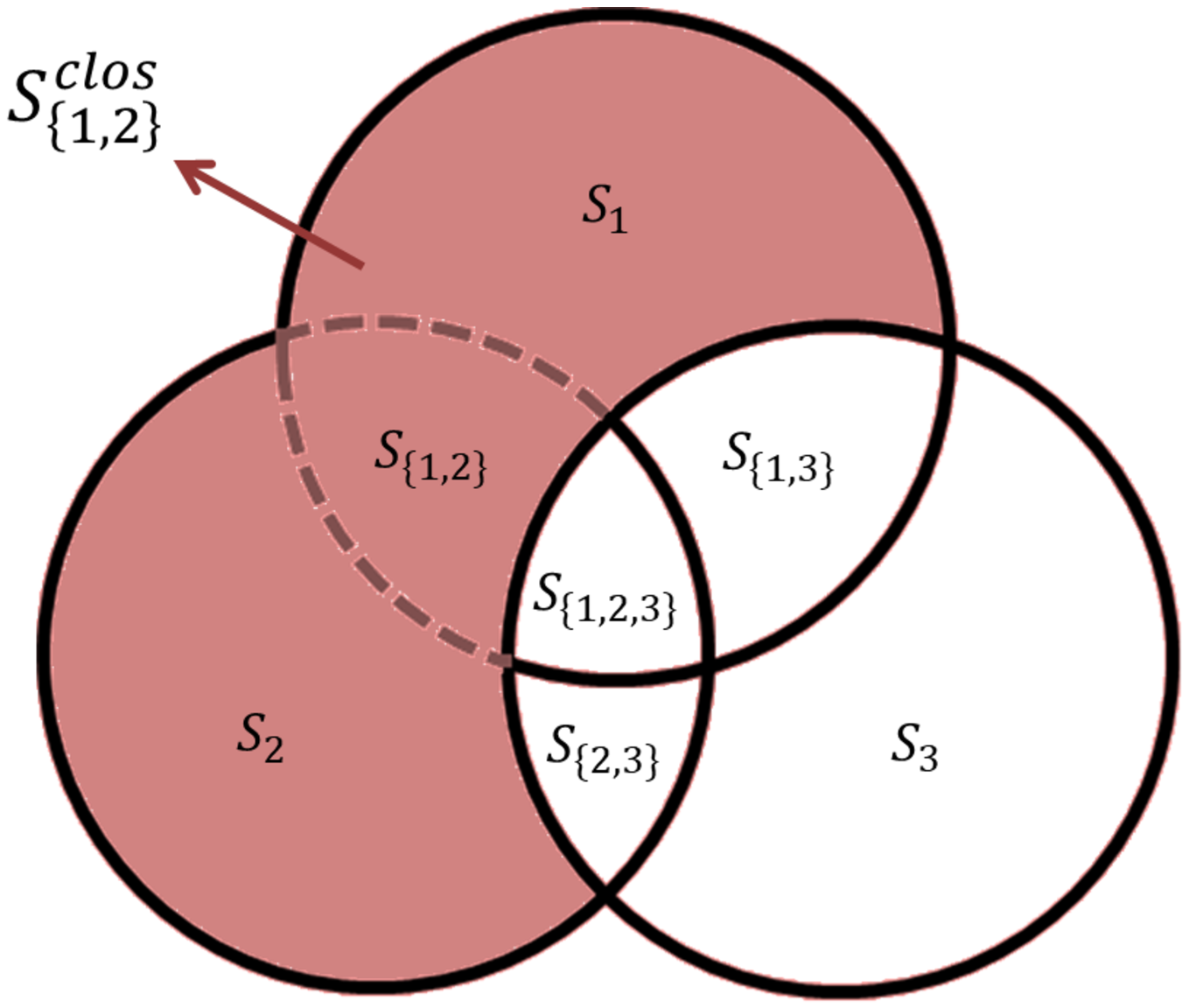}
\caption{Sobol' indices (left) and closed Sobol' indices (right).}\label{fig:sch_sobol}
\end{figure}

Another useful sensitivity index is the closed Sobol' index \cite{Sobol_MMCE_1993} which writes
\begin{equation}
S^{\textrm{clos}}_A = \sum_{B \subset A} S_B = \frac{\mathbb{V}\left(\mathbb{E}\bigl[ G(X) \bigm| X_A\bigr]\right)}{\mathbb{V}[G(X)]}
\end{equation}
In the independent setting, it can be interpreted as the percentage of variability induced by all the variables in a selected subset and their interactions. 
\myfigref{fig:sch_sobol} provides an illustration of the Sobol' indices and the closed Sobol' indices for a model with three inputs. Each Venn diagram represents the variance of the output, with the representation of each of the two Sobol' indices presented above. While this representation is useful in the GSA context, it relies on the assumption of independence between the inputs. 

\section{Axioms of Shapley values}
\label{app:def_shap_val}

Consider a game with $d$ players, and let $\textrm{val}(A) \in \mathbb{R}$ be the cost function quantifying the production value of a coalition (\ie set of players) $A \in \mathcal{P}_d$, under the assumption that $\textrm{val}(\emptyset)=0$. The Shapley value $\phi_j=\phi_j(\textrm{val}), j=1,\dots,d$ attributed to each player can be defined by the following set of axioms:

\begin{enumerate}
\item (Efficiency) $\sum_{j=1}^d \phi_j = \textrm{val}(\{1 \dots, d\})$, meaning that the sum of the allocated values have to be equal to the value produced by the cooperation of all the players.

\item (Symmetry) If $\textrm{val}(A\cup\{ i\}) = \textrm{val}(A\cup \{ j\})$ for all $A \in \mathcal{P}_d$, then $\phi_i = \phi_j$, meaning that if two players allow for the same contribution to every coalition, their attribution should be the same.

\item (Dummy) If $\textrm{val}(A\cup\{ i\}) = \textrm{val}(A)$ for all $A \in \mathcal{P}_d$, then $\phi_i=0$, meaning that if a player does not contribute the the production of resources for all coalition, he should not be attributed any resources.

\item (Additivity) If $\textrm{val}$ and $\textrm{val}'$ have Shapley Values $\phi$ and $\phi'$ respectively, then the game with cost function $\textrm{val} + \textrm{val}'$ has Shapley values $\phi_j + \phi'_j$ for $j \in \{1, \dots, d \}$.
\end{enumerate}

These four axioms guarantee a cooperative allocation of $\textrm{val}(\{1, \dots, d \})$. The unique  attribution method that satisfies these four axioms are the Shapley values \cite{osborne_course_1994}, defined by:
\begin{equation}
\label{eq:shap_val}
\phi_j = \frac{1}{d} \sum_{A \subset -j} {d-1 \choose |A|}^{-1} (\textrm{val}(A\cup\{ j\}) - \textrm{val}(A)), \quad j=1,\dots,d
\end{equation}
where $\{ -j\} = \{ 1, \dots, d \} \backslash j$. One can additionally remark that $\phi_j(\textrm{val})$ is a linear operator, meaning that for some constant $c \in \mathbb{R}$, $\phi_j(c\times \textrm{val}) = c \times \phi_j(\textrm{val})$.

\section{Mathematical proofs}
\label{ap:proofs}

\subsection{Positivity of the $(\ell^1)$-target Shapley effects}\label{ap:pos_proofl1}
Let $A \subseteq \{ 1, \dots, d\} \setminus \{j\}$, for $j \in \{1, \dots, d\}$. In order to show that the $(\ell^1)$-target Shapley effects are positive, one needs to prove that:
\begin{equation}\label{eq:l1TShapley}
\textrm{T-S}_{A \cup \{j\}}^{\ell^1} \geq \textrm{T-S}_A^{\ell^1}.
\end{equation}
In \cite{Cui_Lu_SCTS_2010}, it was shown that the following property holds:
\begin{equation}
\eta_{A \cup \{ j\}} \geq \eta_A
\label{ap:prop_eta}
\end{equation}
with $\eta_A$ being defined in \myeqref{eq:eta}.
From the definition of $\textrm{T-S}^{\ell^1}_A$,
\begin{equation}
\textrm{T-S}^{\ell^1}_A = \frac{2}{\mathbb{E}\Bigl[ \bigl\lvert \mathds{1}_{\mathcal{F}_t}(X) - \mathbb{E}\left[\mathds{1}_{\mathcal{F}_t}(X)\right] \bigr| \Bigr]} \eta_A  ,
\end{equation}
one gets immediately the property \ref{eq:l1TShapley}.


\subsection{Positivity of the $(\ell^2)$-target Shapley effects}\label{ap:pos_proofl2}

Let $X=\left( X_1, \dots, X_d\right) \in \mathbb{R}^d$, be a real-valued random vector admitting a probability measure $P_X$ on the usual real measurable space. Let $L^2(P_X)$ be the functional space such that, for a measurable function $f$, $\displaystyle \LtwoNorm{f} \overset{\textrm{def}}{=} \int_{\mathbb{R}^d} f^2(x) dP_X(x) < + \infty$. Let $G(\cdot) \in L^2$ be the studied numerical model, and denote the random variable $Y=G(X)$ be the model output (or $Y=\mathds{1}_{G(X)>t}(X)$ the TSA variable of interest, without loss of generality). Let $A \subseteq \{1, \dots, d \} \setminus \{j\}$ be the indices of the subset of inputs $X_A$ and $j \in \{1, \dots, d\}$. In order to show that $\textrm{T-Sh}_j \geq 0$, one needs to prove that:
\begin{equation}
\textrm{T-S}_{A \cup \{j\}} -\textrm{T-S}_A \geq 0
\end{equation}
which is equivalent to
\begin{equation}
\VarBigLR{\mathbb{E}\left[ Y | X_A\right]} \leq \VarBigLR{\mathbb{E}\left[ Y | X_{A \cup \{ j\}} \right]}.
\end{equation}

From the Pythagorean theorem, one has:
\begin{equation}
\LtwoNorm{Y} = \LtwoNorm{\mathbb{E}\bigl[ Y\bigm| X_A\bigr]} + \LtwoNorm{Y - \mathbb{E}\bigl[ Y\bigm| X_A\bigr]},
\end{equation}
which is equivalent to
\begin{equation}
\mathbb{E}\bigl[ Y^2 \bigr] = \mathbb{E}\Bigl[ \left(\mathbb{E}\bigl[ Y\bigm| X_A\bigr]\right)^2 \Bigr] + \mathbb{E}\Bigl[ \left( Y - \mathbb{E}\bigl[ Y \bigm| X_A\bigr]\right)^2 \Bigr].
\end{equation}
By removing $\left(\mathbb{E}\bigl[ Y\bigr]\right)^2$ to both sides of the equality, one obtains:
\begin{equation}
\label{ap:eq_var}
\VarBigLR{\mathbb{E}\left[ Y | X_A\right]} = \VarLR{Y} - \LtwoNorm{Y - \mathbb{E}\left[ Y | X_A\right]}.
\end{equation}
By using the formula $\mathbb{E}\bigl[ Y\bigm| X_A\bigr] = \underset{Z \in \sigma\left( X_A\right)}{\textrm{argmin }} \LtwoNorm{Y - Z}$, with $\sigma(X_A)$ being the span of $X_A$, we deduce that $\mathbb{E}\bigl[ Y\bigm| X_A\bigr] \leq \mathbb{E}\bigl[ Y\bigm| X_{A \cup \{j\}}\bigr]$ since $\sigma \bigl(X_A\bigr) \subseteq \sigma \bigl(X_{A \cup \{j\}} \bigr)$. This leads to
\begin{equation}\label{eq:eq_varbis}
\VarLR{Y} - \LtwoNorm{Y - \mathbb{E}\bigl[ Y\bigm| X_A\bigr]} \leq \VarLR{Y} - \LtwoNorm{Y - \mathbb{E}\bigl[ Y\bigm| X_{A\cup \{ j\}}\bigr]}.
\end{equation}
Finally, from \myeqref{ap:eq_var} and \myeqref{eq:eq_varbis}, we obtain
\begin{equation}
\VarBigLR{\mathbb{E}\bigl[ Y\bigm| X_A\bigr]} \leq \VarBigLR{\mathbb{E}\bigl[ Y\bigm| X_{A \cup \{ j\}} \bigr]}
\end{equation}
which concludes the proof.

\section{Minimal R code examples for the estimation methods}
\label{sec:min_code}

\subsection{Monte Carlo sampling estimator}
\begin{lstlisting}[language=R, caption=Minimal R code example for the Monte Carlo estimation.]
#Packages
library(sensitivity)
library(mvtnorm)
library(condMVNorm)

#Model definition
model.linear <- function(X) as.numeric(apply(X,1,sum)>0)

#Parameters
d <- 3
mu <- rep(0,d)
sig <- c(1,1,2)
ro <- 0.9
Cormat <- matrix(c(1,0,0,0,1,ro,0,ro,1),d,d)
Covmat <- ( sig %*% t(sig) ) * Cormat

#Total and marginal simulation function
Xall <- function(n) mvtnorm::rmvnorm(n,mu,Covmat)

#Conditional simulation function
Xset <- function(n, Sj, Sjc, xjc){
  if (is.null(Sjc)){
    if (length(Sj) == 1){ rnorm(n,mu[Sj],sqrt(Covmat[Sj,Sj]))
    }else{ 
    mvtnorm::rmvnorm(n,mu[Sj],Covmat[Sj,Sj])
    }
  }else{ 
  condMVNorm::rcmvnorm(n, 
                    mu, 
                    Covmat, 
                    dependent.ind=Sj, 
                    given.ind=Sjc, 
                    X.given=xjc)
    }
}

#(l2)-target Shapley effects estimation
l2_tse.mc <- shapleyPermEx(model = modlin, 
                            Xall=Xall, 
                            Xset=Xset, 
                            d=d, 
                            Nv=1e4, 
                            No = 1e3, 
                            Ni = 3)
#Plot the results
print(l2_tse.mc)

#(l2)-target Shapley effects estimation with random permutations
l2_tse.mc.randperm<-shapleyPermRand(model = modlin,                        
                            Xall=Xall, 
                            Xset=Xset, 
                            d=d, 
                            Nv=1e4, 
                            No = 1e3, 
                            Ni = 3,
                            m=5)
#Plot the results
plot(l2_tse.mc.randperm)
\end{lstlisting}

\subsection{Nearest-neighbor estimator}
\begin{lstlisting}[language=R, caption=Minimal R code example for the nearest-neighbor estimation.]
#Packages
library(sensitivity)
library(mvtnorm)

#Random sample of inputs-output
X<-rmvnorm(2000, rep(0,3), diag(3))
Y<-rbinom(2000, 1, 0.7)

#(l2)-target Shapley effects estimation
l2_tse.knn<-sobolshap_knn(model=NULL,
              X=X)
tell(l2_tse.knn, Y)

#Plot the results
plot(l2_tse.knn)

#(l2)-target Shapley effects estimation with random permutations
l2_tse.knn.randperm<-sobolshap_knn(model=NULL,
              X=X,
              rand.perm=T,
              n.perm=5)
tell(l2_tse.knn.randperm, Y)

#Plot the results
plot(l2_tse.knn.randperm)
\end{lstlisting}

\rev{
\section{Empirical convergence rate of the estimation scheme}
\setcounter{figure}{0}    
\label{app:empconv}
\subsection{Empirical convergence of the Monte Carlo estimation procedure}
\label{app:empconv_mc}
In order to illustrate the Monte Carlo estimation procedure of the target Shapley effects (see, Section \ref{sec:mc_estim}), the following model is considered: 
\begin{equation}
X = \begin{pmatrix} X_1 \\ X_2 \\ X_3\end{pmatrix}\sim \mathcal{N}\left( \begin{pmatrix} 0 \\ 0 \\ 0 \end{pmatrix}, \begin{pmatrix}1 & 0 & 0 \\ 0 & 1 & 0.6 \\ 0 & 0.6 & 1 \end{pmatrix}\right), \quad Y = \sum_{i=1}^3 X_i
\end{equation}
where the considered TSA variable of interest is $\mathds{1}_{Y>3}(X)$. This represents a failure probability $p_t^Y \simeq 0.071$. The empirical rate of convergence is studied on $100$ repetitions, with respect to several values of $N_v$, with fixed sample sizes $N = 10^5$ and $N_p = 3$. The empirical convergence results are illustrated in \myfigref{fig:emp_conv_mc_fig}.

\begin{figure}[!ht]
\centering
\includegraphics[width=0.8\textwidth]{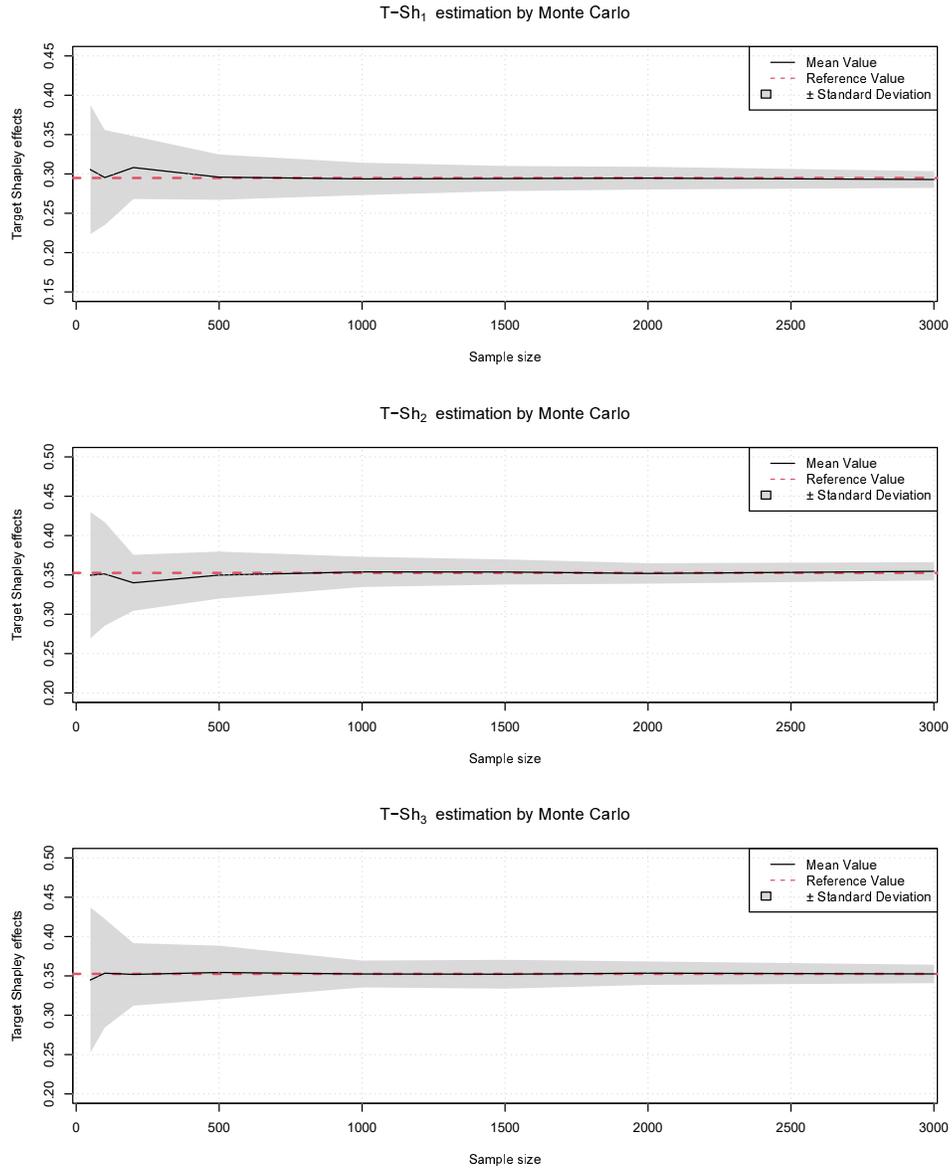}
\caption{Empirical convergence of the Monte Carlo estimation procedure with respect to the sample size.}
\label{fig:emp_conv_mc_fig}
\end{figure}

\subsection{Empirical convergence of the nearest-neighbor estimation procedure}
\label{app:empconv_knn}
The empirical convergence rate of the nearest-neighbor estimator of the target Shapley effects (see, Section \ref{subsec:knn}) is illustrated on the same test-case as in \ref{app:empconv_mc}, with the same TSA variable of interest is $\mathds{1}_{Y>3}(X)$.
The proposed indices have been estimated on $100$ generated samples of $X$, for several sample sizes. The results are presented in \myfigref{fig:emp_conv_knn_fig}.
\begin{figure}[!ht]
\centering
\includegraphics[width=0.8\textwidth]{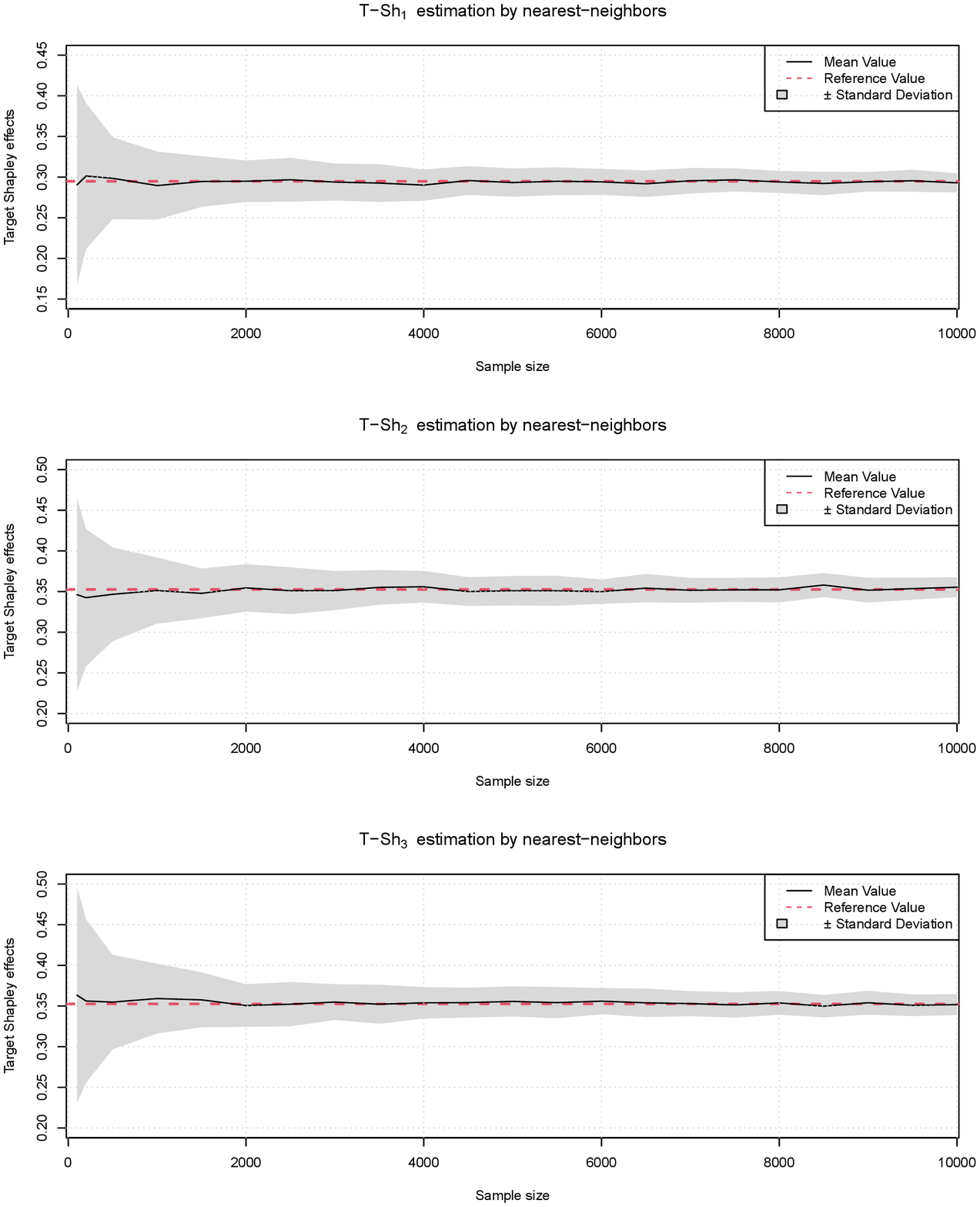}
\caption{Empirical convergence of the nearest-neighbor estimation with respect to the sample size.}
\label{fig:emp_conv_knn_fig}
\end{figure}
}

\end{document}